\documentclass[aop]{imsart}

\RequirePackage{amsthm,amsmath,amsfonts,amssymb}
\RequirePackage[numbers,sort&compress]{natbib}
\RequirePackage[colorlinks,citecolor=blue,urlcolor=blue]{hyperref}
\RequirePackage{graphicx}

\newcommand{\E}{\mathbb E}
\newcommand{\R}{\mathbb R}
\newcommand{\N}{\mathbb N}

\newcommand{\F}{\mathcal F}

\newcommand{\eps}{\varepsilon}

\newcommand{\Ll}{\mathbb{L}}
\newcommand{\ud}{\,\mathrm{d}}

\def\P{{\mathbb P}}

\def\om{\omega}
\def\Om{\Omega}

\newcommand{\Ind}{\mathbf{1}}

\startlocaldefs
\theoremstyle{plain}
\newtheorem{axiom}{Axiom}
\newtheorem{claim}[axiom]{Claim}
\newtheorem{theorem}{Theorem}[section]
\newtheorem{lemma}[theorem]{Lemma}
\newtheorem{corollary}[theorem]{Corollary}
\newtheorem{proposition}[theorem]{Proposition}
\newtheorem{remark}[theorem]{Remark}
\theoremstyle{definition}
\newtheorem{definition}[theorem]{Definition}
\newtheorem*{example}{Example}
\newtheorem*{fact}{Fact}
\endlocaldefs

\begin{document}

\begin{frontmatter}
\title{Hereditary Hsu-Robbins-Erd\H{o}s    \\ 
Law of Large Numbers}

\runtitle{Hereditary Hsu-Robbins-Erd\H{o}s LLN}
 \thankstext{T1}{~To  the memory of Herbert Ellis   Robbins  (1915--2001).}

\begin{aug}

\author[A]{\fnms{Istv\'an}~\snm{ Berkes} 
\ead[label=e1]{berkes.istvan@renyi.hu}},
\author[B]{\fnms{Ioannis}~\snm{Karatzas}\ead[label=e2]{ik1@columbia.edu}\orcid{0000-0000-0000-0000}}
\and
\author[C]{\fnms{Walter}~\snm{Schachermayer}\ead[label=e3]{walter.schachermayer@univie.ac.at}}
%%%%%%%%%%%%%%%%%%%%%%%%%%%%%%%%%%%%%%%%%%%%%%
%% Addresses                                %%
%%%%%%%%%%%%%%%%%%%%%%%%%%%%%%%%%%%%%%%%%%%%%%
\address[A]{~Alfred R\'enyi Institute of Mathematics, Re\'altanoda utca 13-15, 1053 Budapest, Hungary  \textbf{%[Additional affiliations should be put in the Acknowledgments section]
}\printead[presep={ ,\ }]{e1}}

\address[B]{~Departments of Mathematics and Statistics,  Columbia University, New York, NY 10027 \textbf{%[Additional affiliations should be put in the Acknowledgments section]
}\printead[presep={,\ }]{e2}}

\address[C]{~Faculty of Mathematics, University of Vienna, Oskar-Morgenstern-Platz 1, 1090 Vienna, Austria \textbf{%[Additional affiliations should be put in the Acknowledgments section]
}\printead[presep={,\ }]{e3}}
\end{aug}

\medskip
\begin{abstract}
We show that every sequence $f_1, f_2, \cdots$ of real-valued random variables with $\,\sup_{n \in \N} \E (f_n^2) < \infty\,$  contains a  subsequence $f_{k_1}, f_{k_2}, \cdots$  converging  in \textsc{Ces\`aro} mean to some $\,f_\infty \in \mathbb{L}^2$ {\it completely,}  to wit, 
$$
\sum_{N \in \N} \, \P \bigg( \bigg| \frac{1}{N} \sum_{n=1}^N   f_{k_n} - f_\infty \bigg| > \eps   \bigg)< \infty\,, \quad \forall ~ \eps > 0\,;
$$
and {\it hereditarily,}  i.e.,   along   all further subsequences and    permutations. We also   identify a condition, slightly weaker than boundedness in  $\, \mathbb{L}^2,$ which turns out to be not only sufficient      for the above hereditary complete convergence in \textsc{Ces\`aro} mean, but    necessary as well.
\end{abstract}



\end{frontmatter}
%%%%%%%%%%%%%%%%%%%%%%%%%%%%%%%%%%%%%%%%%%%%%%
%% Please use \tableofcontents for articles %%
%% with 50 pages and more                   %%
%%%%%%%%%%%%%%%%%%%%%%%%%%%%%%%%%%%%%%%%%%%%%%
%\tableofcontents

\section{Introduction}

The strong law of large numbers (SLLN; \textsc{Kolmogorov} \cite{K2},\,\cite{KAN}; also \cite{Du}) is one of the pillars of   probability theory. For a sequence of real-valued and integrable functions $f_1, f_2, \cdots \,,$ defined on a probability space $(\Omega, \mathcal{F}, \mathbb{P})$ and   independent  with common distribution $   \mu \, $, it  states that the sample,   or "\textsc{Ces\`aro}", averages
\begin{equation}
\label{1.0}
  \frac{1}{N} \sum_{n=1}^N f_{n}\,\,, \quad N \in \N\,,
\end{equation}
converge  $\mathbb{P}-$a.e.\,to the ensemble average  
$\,\E (f_1) = \int_\R \,x\,     \mu   (\mathrm{d} x),$ as $\,N \to \infty\,.$

In 1947, \textsc{Hsu \& Robbins} \cite{HR} obtained a remarkable strengthening of this result. In the same setting as above,  but now under the   square-integrability condition 
\begin{equation}
\label{1.1}
\,\E \big(f_1^2\big) = \int_\R \,x^2\,      \mu    (\mathrm{d} x) < \infty 
\end{equation}
and with $\E \big( f_1 \big) =0$, they established the stronger (so-called "complete") convergence
\begin{equation}
\label{1.2}
\sum_{N \in \N} \, \P \left( \bigg| \frac{1}{N}\sum_{n=1}^N   f_{ n}  \bigg| > \eps   \right)< \infty\,, \qquad \forall ~ \eps > 0 \,.
\end{equation} 
Then in 1949/50, \textsc{Erd\H{o}s} (\cite{Erd},\,\cite{Erd1}) gave a very concise proof of the \textsc{Hsu-Robbins} theorem; and showed  additionally that the square-integrability condition \eqref{1.1} is not only sufficient  for the validity of \eqref{1.2} but also {\it necessary,} as   had been conjectured in \cite{HR}. 

A useful juxtaposition of these results comes about, when one considers the sojourn times
\begin{equation}
\label{1.3}
T_\eps := \sum_{N \in \N} \mathbf{ 1}_{ \big\{ \big| \sum_{n=1}^N    f_{ n}      \big|> \eps N \big\} }  \,, \quad  \eps > 0
\end{equation} 
  the sequence of averages in \eqref{1.0} spends outside $\eps-$neighborhoods of the ensemble average $\,\E (f_1  )=0\,$. For a sequence of independent and equi-distributed $f_1, f_2, \cdots \,$, the SLLN amounts to the statement
\begin{equation}
\label{1.4}
\E \big( \big| f_1 \big| \big ) < \infty \, \,\,\Longleftrightarrow \,\,\, \P \big(T_\eps < \infty \big) =1\,,~~ \forall ~ \eps >0\,;
\end{equation} 
the \textsc{Hsu-Robbins} theorem to the statement 
\begin{equation}
\label{1.5}
\E \big(   f_1^2 \big ) < \infty \, \,\,\Longrightarrow \,\,\, \E \big(T_\eps \big) < \infty   \,,~~ \forall ~ \eps >0\,;
\end{equation}
and   \textsc{Erd\H{o}s}'s result  to the validity   of the reverse implication in \eqref{1.5}.   The connection between the  condition $\E \big(   f_1^2 \big ) < \infty $ of \eqref{1.1}, and complete convergence, was further sharpened by \textsc{Heyde} \cite{H} (see also  \cite{Pr})   who showed that, under \eqref{1.1}, the variance   $\sigma^2 :=  \E \big(   f_1^2 \big )  $ admits   the representation
 \begin{equation} 
 \label{1.6}
\sigma^2\,=\, \lim_{\eps \downarrow 0} \,\left( \eps^2 \sum_{N \in \N} \, \P \Big( \Big|  \sum_{n=1}^N   f_{ n}  \Big| > \eps \,N  \Big) \right) = \, 
\lim_{\eps \downarrow 0} \,\Big( \eps^2 \cdot \E \big(T_\eps \big) \Big)  \,.
\end{equation}

  \smallskip
  \noindent
 
Another classical result along a similar, yet somewhat distinct, strand of inquiry,          is that of \textsc{Koml\'os} \cite{K}. Also known for a long time (60 years already) but always very striking,  it asserts that such "ergodicity" (stabilization via  averaging)  as manifest in the SLLN,  occurs within {\it any} sequence $f_1, f_2, \cdots \,$ of measurable, real-valued  functions  that satisfy  only the  boundedness-in-$\mathbb{L}^1$ condition   
\begin{equation} 
\label{1.7} 
  \sup_{n \in \N} \, \E \big( | f_n| \big) < \infty\,.
  \end{equation}
   More precisely, under \eqref{1.7} {\it   there exist  an integrable function $ f_\infty$ and a subsequence $ f_{k_1}, f_{k_2} , \cdots   $ with 
   \begin{equation} 
\label{1.8} 
\lim_{N \to \infty}   \frac{1}{N} \sum_{n=1}^N f_{k_n} = f_\infty\,, ~~~\text{$\,\mathbb{P}-$a.e.};
  \end{equation}
  and the same holds   "hereditarily", i.e., along every further subsequence  of   $ f_{k_1}, f_{k_2} , \cdots   $} (this  hereditary  aspect is   immediate, when the $f_1, f_2, \cdots \,$ are independent and equi-distributed). Conditions slightly weaker than \eqref{1.7} have been  identified recently in \cite{BKS} as not only sufficient, but also necessary, for the validity of the hereditary SLLN  in \eqref{1.8} (cf.\,\cite{L}).   
   
   The proof     in \cite{K} is an  early application     in probability limit theory   of martingale techniques; these were then used in \cite{Ch1}--\cite{Ch3}, \cite{G} to derive similar hereditary results for the central limit theorem and for the law of the iterated logarithm.

   %%%%%%%%%%%%%
   \subsection{Preview}
   \label{sec1.1}
   %%%%%%%%%%%%%
   Our result, a special case of which is stated succinctly as the first sentence of the abstract and as Corollary  \ref{H-HRE-Cor}, appears       in Theorem \ref{H-HRE}. In the spirit of \textsc{Koml\'os} \cite{K}, it provides a hereditary version   for  the \textsc{Hsu-Robbins-Erd\H{o}s} theorem. This consists of a {\it sufficiency part} (i), corresponding to \cite{HR} and proved in sections \ref{sec6},\,\ref{sec5}; and of a {\it necessity part} (ii), corresponding to \cite{Erd} and proved in section \ref{sec3d}.   
    
    We establish first  an   auxiliary sufficiency result, Proposition  \ref{H-HRE-1},  under an  additional requirement  ($\star$)   which  posits the existence of a subsequence  $ f_{k_1}, f_{k_2} , \cdots   $  whose squares converge  weakly in $\Ll^1$ to some    bounded function.     This auxiliary  result is proved    in section \ref{sec4} via     uniform integrability arguments based on   the KPR   Lemma \ref{KPR} (\cite{KP},\,\cite{Ro}; see also\,\cite{CM},\,\cite{KK}),      perturbation methodologies, and martingale techniques  as in  \cite{K} (cf.\,\cite{Ch},\,\,pp.\,137-141). These  are reviewed  in an Appendix (section \ref{sec7}),  and  are used to effect reductions to   simple martingale differences. A uniform version of the  result    \eqref{1.6},      important in the present context, is  established in  section \ref{sec8}.  
     
    We then proceed to  the proof of Theorem \ref{H-HRE}, using Proposition \ref{H-HRE-1} as a stepping stone. Here, it turns out to be of crucial importance, indeed indispensable,  to go {\it beyond martingale techniques, and to deploy arguments based on exchangeability}.

    %%%%%%%%%%%%%%%%%
     \subsection{Additional Aspects}
\label{sec1.2}
%%%%%%%%%%%%%%%%%%

The statement of Corollary  \ref{H-HRE-Cor}  formulates, in the context of the \textsc{Hsu-Robbins} result,  the heuristic "general principle of subsequences" which appears on the first pages of \textsc{Chatterji} \cite{Ch3} and of \textsc{Berkes-P\'eter} \cite{BerPet} (cf.\,\cite{A1},\,\cite{Ch2}). Its proof proceeds by approximating appropriate subsequences  of $f_1, f_2, \cdots \,$ by sequences      {\it strongly exchangeable at infinity;} pioneered by   \textsc{Aldous} (\cite{A},\,\cite{A1}),  this approach was refined in \textsc{Berkes-P\'eter}\,\cite{BerPet} and is adapted  to our $\Ll^2$ setting here in  somewhat simplified form (section \ref{sec6}).      Exchangeability methods have been deployed (cf.\,\cite{A},\,\cite{A1},\,\cite{BerPet})    to establish instances of the  subsequence  principle for "almost-sure" and for distributional results; but   not for complete convergence   as  done       here, or for convergence  in probability as done  in \cite{BKS}, \cite{KS}.  

An additional     hereditary feature   of the \textsc{Koml\'os}    result  was established in 
 \textsc{Berkes-Tichy}   \cite{B},\,\cite{BerTic}: not only does every subsequence of $\, f_{k_1}, f_{k_2}, \cdots\,$  converge  a.e.\,in \textsc{Ces\`aro} mean  to   $f_\infty  \,$, but so do all further   permutations.  This holds also in our setting,  {\it and all our hereditary results are phrased in terms of both subsequences and permutations}.

Yet another    noteworthy feature occurs here.   \textsc{Beno\^ist-Quint} \cite{BQ} construct   a bounded-in-$\Ll^2$   martingale-difference sequence  with $\,\sum_{N \in \N} \, \mathbb{P}  ( \big|  \sum_{n=1}^N   f_{n}   \big| > \eps\, N   )= \infty\,$ for all  $\, \eps > 0\,$  (but   containing, in accordance with Corollary \ref{H-HRE-Cor} below, a   subsequence  $f_{k_1}, f_{k_2}, \cdots$ along which, and along whose every further subsequence,     $\,\sum_{N \in \N} \, \mathbb{P}  ( \big|  \sum_{n=1}^N   f_{k_n}   \big| > \eps\, N   )< \infty\,$ holds for all $\,     \eps > 0\,$). Thus, for proving   Theorem \ref{H-HRE}\,(i),   or even its Corollary \ref{H-HRE-Cor}, {\it approximation via martingale differences alone  leads to  a dead end;} rather, one needs  to approximate a suitable subsequence by an  {\it exchangeable} sequence.  

That such a situation should arise, had   
been predicted by \textsc{Aldous} in \cite{A},\,\cite{A1}; but the argument  was  supported there   only via (very) artificial examples. Our  Theorem \ref{H-HRE}     seems to provide  the first "concrete" instance, of an important  limiting result   in Probability Theory    valid  for independent, equi-distributed functions,   whose  hereditary extension  {\bf requires}     methods based on exchangeability.

\smallskip

%%%%%%%%%%%%%%%%%%%%%%%
\section{A Hereditary \textsc{Hsu-Robbins-Erd\H{o}s}  Law of Large Numbers}
\label{sec2}
%%%%%%%%%%%%%%%%%%%%%%%

\begin{definition} {\it The HRE Property.} 
\label{D-HRE}
Consider a sequence of real-valued, measurable functions $\,  f_1, f_2, \cdots \,$ on a probability space $(\Om, \mathcal{F}, \mathbb{P})$.     We say that this sequence  {\it  satisfies the HRE} (\textsc{Hsu-Robbins-Erd\H{o}s}) {\it property} for some real-valued,  measurable function   $f_\infty \, ,$ if it contains   a subsequence $ f_{k_1}, f_{k_2}, \cdots $   converging in \textsc{Ces\`aro} mean to  this $f_\infty  $ completely, i.e., 
\begin{equation} 
\label{2.4}
\sum_{N \in \N} \, \P \Big( \Big|   \sum_{n=1}^N   f_{{k_n}} - N\,f_\infty   \Big| > \eps \,N  \Big)< \infty\,, \qquad \forall ~ \eps > 0\,,
\end{equation} 
and  hereditarily,  i.e., also along all   subsequences   of $ f_{k_1}, f_{k_2}, \cdots $ and all their permutations. 
\end{definition}

 We    establish in this paper the following hereditary version of the \textsc{Hsu-Robbins-Erd\H{o}s}   Law of Large Numbers (\cite{HR}, as well as \cite{Erd},\,\cite{Erd1}).
   
\begin{theorem} 
[\textsc{Hereditary Hsu-Robbins-Erd\H{o}s  LLN}]
\label{H-HRE}
On a probability space $(\Omega, \mathcal{F}, \mathbb{P})$, consider a sequence of real-valued, measurable  functions  $\,f_1, f_2, \cdots\,.$

 \smallskip
 \noindent
 {\bf (i)}  Suppose that    for some sequence of  sets     $  \,A_1\,, \, A_2 \,, \cdots $ in $\,  {\cal F}$ with $\, \lim_{n \to \infty }  \P (A_n) =1  \,,$   \\ {\bf .}       the sequence  $\,f_{1} \mathbf{ 1}_{A_{1}}, \,f_{ 2}\, \mathbf{ 1}_{A_{2}}, \cdots $  is  bounded in $\Ll^2\,;$      while \\ {\bf .}    the sequence $\,f_{1} \mathbf{ 1}_{A_{1}^c}, \,f_{ 2}\, \mathbf{ 1}_{A_{2}^c}, \cdots $ converges to  zero in  $\,\Ll^1\,.$

\noindent
 There exists  then   a    subsequence $ f_{k_1}, f_{k_2}, \cdots $ of $ \,f_1, f_2, \cdots \, $ with   the HRE property for  some   $f_\infty  \in \Ll^2\,.$ 
   
 \smallskip
 \noindent
 {\bf  (ii)}  Conversely, suppose     $\, f_1, f_2, \cdots \,$   
 has the HRE      property  for some   $f_\infty \in \Ll^2\,.$       There exist  then  a  sequence         $  \,A_1\,, \, A_2 \,, \cdots $ of  sets in $\,  {\cal F}$    with $\, \lim_{n \to \infty }    \P (A_n) =1  \,,$ and a    subsequence $\, f_{k_1} ,f_{k_2},\cdots$  of $\, f_{ 1} ,f_{ 2},\cdots\,,$  such that   \\ {\bf .}   the sequence $\, f_{k_1} \mathbf{ 1}_{A_{k_1}}, \, f_{k_2} \mathbf{ 1}_{A_{k_2}}, \cdots $    is bounded in $\, \Ll^2,$  and  \\  {\bf .}    the sequence   $\,f_{ k_1} \mathbf{ 1}_{A_{k_1}^c},\, f_{ k_2} \mathbf{ 1}_{A_{k_2}^c}, \cdots \,$ converges to zero  in $\,\Ll^1.$
\end{theorem} 

If all its  sets $\, A_1 \,,\, A_2 \,, \cdots\,$    
satisfy $\, \P (A_n) =1\,,$      Theorem \ref{H-HRE}\,(i)  amounts to the   statement that follows. This corresponds to     the first sentence of  the paper's abstract, and vindicates in the present context the heuristic   principle of subsequences enunciated  in   \cite{Ch2},\,\cite{Ch3} (cf.\,\cite{BerPet}). 

\begin{corollary} 
\label{H-HRE-Cor}
Every  sequence of real-valued, measurable functions $f_1, f_2, \cdots \,$  which is bounded in $\,\mathbb{L}^2\,$, i.e., satisfies
 \begin{equation} 
 \label{2.1}
 \sup_{n \in \N} \,\E \big(   f_n^2\big) < \infty\,,
 \end{equation}
  contains a subsequence   $ ~f_{k_1}, f_{k_2} , $ $\cdots  $  with the HRE property for some $f_\infty \in \Ll^2$.   
\end{corollary}

In section \ref{sec4} we shall establish,      using   martingale methods,  the preliminary  result of    Proposition \ref{H-HRE-1} below   with  its additional assumption $(\star)\,$; then deploy this result in sections \ref{sec5},\,\ref{sec3d}  as a crucial  stepping stone that allows us eventually   to ascend to the generality of Theorem \ref{H-HRE}.

\begin{proposition} 
\label{H-HRE-1}
Suppose that the    sequence of real-valued, measurable functions $f_1, f_2, \cdots \,$ is bounded in $\,\mathbb{L}^2\,$ as in \eqref{2.1}, and satisfies the following condition:

\noindent
 $(\star)$    There exists a subsequence $ f_{k_1}, f_{k_2} , \cdots   $ of $f_1, f_2, \cdots \,,$  whose squares $ f_{k_1}^2, f_{k_2}^2 , \cdots   $ converge weakly in $\mathbb{L}^1$ to a function  $\,  \boldsymbol{\eta}  \in \Ll^\infty ,$ namely,  
\begin{equation} 
\label{2.2}
\lim_{n \to \infty} \E \big( f_{k_n}^2 \cdot \xi \big) = \E \big(   \boldsymbol{\eta}  \cdot \xi \big)\,, \quad \forall ~~ \xi \in \Ll^\infty\,.
\end{equation}
 There exist then  a real-valued function $
f_\infty  \in \Ll^2$ and a suitable (further, relabelled)  subsequence   $ f_{k_1}, f_{k_2} , \cdots  $ of $\,f_1, f_2, \cdots $   with the HRE property for this $f_\infty\,$. 
 \end{proposition}

\begin{remark}
{\it Stable Convergence.} 
\label{rem2.4a}
{\rm
The  functions $ \,f_\infty\,, \,  \boldsymbol{\eta} \,$        play the r\^oles of randomized limiting first and second moments  for the   subsequence $f_{k_1}, f_{k_2}, \cdots \,$; they are measurable with respect to the tail $\sigma-$algebra   
\begin{equation} 
\label{2.3tail}
{\cal T}  := \bigcap_{n \in \N} {\cal T}_n\,,\qquad {\cal T}_n := \sigma \big( f_{k_n}, f_{k_{n+1}}, \cdots \big).
\end{equation}

Let us elaborate. We recall  the   space $\Ll^0$ of real valued,  measurable functions, endowed with convergence in probability.    As shown in \S\,\ref{sec6.1.1} here, every   sequence $\,f_1, f_2, \cdots\,$ with the HRE property contains a subsequence $\, f_{k_1} ,f_{k_2},\cdots$  which  is  {\it bounded in} $\Ll^0$ (bounded in probability, tight): namely, with $\, \sup_{n \in \N} \,\P \big( \big| f_{k_n} \big| > \lambda \big) \rightarrow 0\,$ as   $\, \lambda  \to \infty\,.$  

Thus,  for   Theorem \ref{H-HRE} we may assume    $f_1, f_2, \cdots \,$ to contain a    subsequence bounded in $\Ll^0$. Now, from a long line of inquiry initiated by \textsc{R\'enyi} in \cite{Ren} (cf.\,\cite{AE};\,\cite{BC};\,\cite{BerRos},\,Theorem 2.2),   such  a sequence   contains also a (relabelled) "determining"   subsequence $ f_{k_1} , f_{k_2}  , \cdots   ,$ along which    the   {\it stable convergence} (extended  \textsc{Helly-Bray} lemma)    
\begin{equation} 
\label{2.7}
\lim_{n \to \infty} \P \big( f_{k_n}  \le x, B \big) =:   \,\mathbf{Q}  (x, B) = \int_B   \mathbf{H}  (x, \omega) \, \P (\mathrm{d} \omega) \,, \qquad \forall ~ ~B \in \F  
\end{equation}
 holds  at every point $x  $ of a countable,   dense  set $ \, \mathbf{D}  \subset \R\,.$   
 
 Here, for each   $x \in   \mathbf{D} $,  the map  $B \mapsto \mathbf{Q}  (x, B) $  is a measure on $\F$,   absolutely continuous with respect to $\P\,$;  and    $\,\omega \mapsto  \mathbf{H}  (x, \omega)$   a version of the \textsc{Radon-Nikod\'ym} derivative $ \mathrm{d}  \mathbf{Q}  (x, \cdot) /  \mathrm{d}  \mathbb{P}\,,$ measurable with respect to the tail $\sigma-$algebra  $\,  {\cal T}  \, $ in \eqref{2.3tail}. 
   Whereas, for $\P-$a.e.\,$\,\omega \in \Omega\,$, the map    $\,x \mapsto  \mathbf{H}   (x, \omega) $ is a probability distribution function on the real line    (cf.\,\cite{A},\,\,Lemma 2),      called  {\it limit random distribution function of the determining subsequence} $ \,f_{k_1} , f_{k_2}  , \cdots   $. We denote its first   (when it exists) and   second moments, respectively,    by 
\begin{equation} 
\label{2.8}
f_\infty (\omega) = \int_\R x \, \mathrm{d}   \mathbf{H}   (x, \omega)      \,, \quad
 \boldsymbol{\eta}  (\omega) = \int_\R  x^2 \,   \mathrm{d}   \mathbf{ H}   (x, \omega)  \,, \qquad \omega \in \Omega \,.
\end{equation}
The function  $\,f_\infty: \Omega \to \R\,$    is well-defined and integrable under the condition \eqref{1.7} of the \textsc{Koml\'os} theorem, and then    \eqref{1.8} holds; whereas, $\,  \boldsymbol{\eta}  : \Omega \to [0, \infty]\,$ is the function   in Proposition \ref{H-HRE-1}.    These same  $\,f_\infty (\omega)\,, \, \boldsymbol{\eta}  (\omega)\,, $ as well as the distribution function $\,   \mathbf{ H }  (\cdot\,, \omega)\,  $ and the   measure $\,   \boldsymbol{\mu}  (\omega)\,$  it induces  on the \textsc{Borel} sets of the real line, are generated also by the sequence   $\,f_{1} \mathbf{ 1}_{A_{1}}, \,f_{ 2}\, \mathbf{ 1}_{A_{2}}, \cdots $.
}
 \end{remark}

\begin{remark} 
\label{rem2.5a} 
{\rm
 When  the distribution function     $x \mapsto   \mathbf{H}  (x )$  in \eqref{2.7} happens   to be independent of $\omega$, the stable convergence $\, \lim_{n \to \infty} \P \big( f_{k_n}  \le x, B \big)  =     \mathbf{H}  (x   ) \,  \P (B) \,, ~\, \forall   \, \,B \in \F \,$  in \eqref{2.7} is called   {\it mixing} (e.g.,\,\cite{AE},\,\cite{Ren}). The quantities   of   \eqref{2.8} are  then  real  constants, and the condition  $\,  \boldsymbol{\eta}   \in \Ll^\infty\,$ of Proposition \ref{H-HRE-1} is satisfied trivially.   
}
\end{remark}

\begin{remark} 
\label{rem2.5} 
{\rm
The conditions of Proposition \ref{H-HRE-1}   hold in  the original \textsc{Hsu-Robbins-Erd\H{o}s}  Law of Large Numbers, when the  $f_1, f_2, \cdots \,$ are independent, equi-distributed  and square-integrable with $  \E (f_1 )=0\,$; for then we can take $\, \P (A_n) =1\,,\, \, \forall $  $n \in \N\,$ in part {\bf (i)} of Theorem \ref{H-HRE}  and   obtain \eqref{2.1} trivially, \eqref{2.2}   with  $   \boldsymbol{\eta}  = \E (f_1^2)$, and  \eqref{2.4}   with $f_\infty = 0.$ 
}
\end{remark}

%%%%%%%%%%%%%%%%%%%%%%
\subsection{Martingale-Difference Sequences}
\label{sec3md}
%%%%%%%%%%%%%%%%%%%%%%%

Suppose       $f_1, f_2, \cdots \,$ is a martingale-difference sequence  with respect to its own filtration.    If   it   is also bounded in $\Ll^{p}$ for some $p >2\,$ (this requirement is stronger than    \eqref{2.6}),   Theorem 3.6 in \textsc{Lesigne\,-Voln\'y}  \cite{LV}     shows    that the   $f_1, f_2, \cdots \,$ converge completely in \textsc{Ces\`aro} mean  to zero: i.e., \eqref{2.4} holds   with $f_\infty \equiv 0$.

Our focus here, is on   the HRE case $p=2\,$.   In a significant  development, mentioned already in subsection \ref{sec1.2},  \textsc{Beno\^ist-Quint} construct  (cf.\,Remark 2.3  in\,\cite{BQ}) a  martingale-difference sequence   $f_1, f_2, \cdots \,, $     bounded  in $\Ll^2$, \, {\it for which the complete convergence of \eqref{2.4}    with $f_\infty \equiv 0\,$  fails;}  they   point out that the functions   $f_1, f_2, \cdots \, $ may even be   independent. It is important to stress that, in     references \cite{BQ} and \cite{LV},  the results refer {\it to the martingale difference sequences themselves, without any passage to subsequences.}

 %%%%%%%%%%%%%%%%%%%%%%
\subsection{Exchangeable Sequences}
\label{sec3c}
%%%%%%%%%%%%%%%%%%%%%%%

Suppose   the $f_1, f_2, \cdots$ are   exchangeable, i.e., their finite-dimensional marginal distributions are invariant under permutations of indices, with $\E(|f_1|)<\infty$. The celebrated    \textsc{de\,Finetti} theorem (e.g.,\,\cite{CT},\,\cite{King})     
provides then a random probability distribution function  $x \mapsto  \mathbf{ H}  (x, \omega)$, measurable with respect to the tail $\sigma-$algebra $\, {\cal T}\, $   and satisfying    $ \, \int_\R  |x| \, \mathrm{d}   \mathbf{ H} (x, \omega)$ $ < \infty\,$ for $\,\P-$a.e. $\omega \in \Omega\,,$ with the following  property:  {\it "\,Given  $\,{\cal T} ,$   the    $f_1, f_2, \cdots$ are independent,  with $\,\P \big[ f_1 \le x \, \big|\, {\cal T}  \big] (\omega) =   \mathbf{ H}   \big( x;\omega \big),$ $\,x \in  \R \,$ as their common distribution function,   \,  for $\,\P-$a.e.\,\,$\omega \in \Omega$."} \, The first and second moments    
\begin{equation} 
\label{2.8_too}
~f_\infty (\omega)=\E \big[  f_1    \big| {\cal T}   \big] (\omega) = \int_\R x \,    \mathrm{d}  \mathbf{ H}   \big(   x ; \omega \big), \quad  \boldsymbol{\eta}  (\omega) =\E \big[   f_1^2     \big|  {\cal T}  \big] (\omega)= \int_\R x^2 \,  \mathrm{d}  \mathbf{ H}   \big(   x ; \omega \big),
\end{equation} 
of this random distribution,     provide here the functions in \eqref{2.8}, and we have  $\E (   \boldsymbol{\eta} ) = \E (f_1^2) $.    For   exchangeable  $f_1, f_2, \cdots$,  it follows from Proposition \ref{prop8.0} that    $ \E (f_1^2) < \infty  $ is equivalent to  
 \begin{equation} 
\label{2.4b}
\sum_{N \in \N} \, \P \left( \bigg| \frac{1}{N} \sum_{n=1}^N   f_{ n} - f_\infty   \bigg| > \eps   \right)< \infty\,, \qquad \forall ~ \eps > 0\,,
\end{equation}
and to the following  strengthening, inspired by \textsc{Heyde}'s identity \eqref{1.6},   of \eqref{2.4b}: 
 \begin{equation} 
\label{2.9}
 \limsup_{  \eps \downarrow 0} \left[ \,  \eps^{2} \sum_{N \in \N} \, \P \left( \bigg|  \frac{1}{N}  \sum_{n=1}^N   f_{ n}   - f_\infty   \bigg| > \eps      \right)  \, \right] \,<\, \infty\,.
 \end{equation} 
  
 \begin{proposition} 
\label{H-HRE-Three}
 Let $\,f_1, f_2, \cdots\,$ be an  
 exchangeable   sequence of real-valued, integrable functions, and let the functions $f_\infty\,,\,   \boldsymbol{\eta} \,$ be as in \eqref{2.8_too}.   
 
  Then  \,\eqref{2.4b}  $\Longleftrightarrow$ \eqref{2.9}  $\Longleftrightarrow \, f_1 \in \Ll^2  \, \Longleftrightarrow\, \E (   \boldsymbol{\eta} ) <\infty.$
 \end{proposition}

 %%%%%%%%%%%%%%%%%%%%%%
\subsection{Uniform Integrability Considerations}
\label{rem2.4}
%%%%%%%%%%%%%%%%%%%%%%%

For  establishing the HRE property \eqref{2.4}   in the context of    Theorem \ref{H-HRE}\,(i),      we may assume     that  
\begin{equation} 
\label{2.6}
 \text{the sequence $\,f_1^2 \,\mathbf{ 1}_{A_{1}}, f_2^2 \, \mathbf{ 1}_{A_{2}}, \cdots \, $ is uniformly integrable.}
  \end{equation}

To justify this   claim,   let us   place ourselves in the setting of Theorem \ref{H-HRE}\,(i).    The sequence   $\,  ( f_n^2 \cdot \Ind_{A_n }  )_{n \in \N}\,$ is bounded in $\Ll^1,$ so    Lemma \ref{KPR} provides a sequence $B_1, B_2, \cdots$ of disjoint sets in ${\cal F}$ such that, after passing to a    subsequence,  $\, \big( f_n^2 \cdot \Ind_{A_n \setminus B_n } \big)_{n \in \N}\,$ is uniformly integrable.  On the other hand, the functions      $\,  h_n:=  f_n  \cdot \Ind_{A_n \cap B_n }\,,\,n \in \N \,$    are bounded in $\Ll^2 $ and supported on disjoint sets,    hence tending to zero in the norm of $\Ll^1$. Invoking Proposition \ref{lem3.7}\,(i) once again we deduce, possibly after passing to a subsequence, that this sequence satisfies the HRE property with $f_\infty \equiv 0\,$. It suffices,  therefore, to establish  this HRE property for the sequence of functions  $\, \big( f_n  \cdot \Ind_{A_n \setminus B_n } \big)_{n \in \N}\,,$   whose squares are uniformly integrable.

 %%%%%%%%%%%%%%%%%%%%%%
\subsection{Ramifications}
\label{sec2bk}
%%%%%%%%%%%%%%%%%%%%%%%

Suppose $f, f_1, f_2, \cdots$ are independent real-valued measurable functions with the same distribution, and let $\, 0 < r < 2\,$, $ p \ge r$. \textsc{Baum \& Katz} \cite{BK} proved that 
\begin{equation} 
\label{2.6z}
\sum_{N \in \N}\,N^{(p/r)-2}\, \,\P \left( \bigg|    \sum_{n=1}^N   f_{ n}      \bigg| > \eps  \,N^{1/r}     \right)\, < \, \infty\,, ~~~ \forall ~ \eps >0
  \end{equation}
holds if, and only if, $\E \big(  | f  |^p\big) < \infty;$ and, if $p \ge 1$, then $  \E  ( f  ) =0 .$ The  \textsc{Hsu-Robbins-Erd\H{o}s} result corresponds to $\,p=2\,, \,r=1.$

In a similar vein, \textsc{Stoica} \cite{St} states a hereditary version of \eqref{2.6z} for any $\Ll^p-$bounded sequence $  f_1, f_2, \cdots$   with $0<r\le p <2\,.$ The argument  in  \cite{St} is fairly direct (it needs also $\, \E  ( f_n  ) =0\, $ when $p \ge 1$)\,;    but does not cover the HRE case investigated here. 
  
  This underscores the fact that, in the \textsc{Baum-Katz} theorem, there   are  substantial differences among the cases $p > 2 \,$ already discussed in subsection \ref{sec3md},    the considerably simpler  $p < 2$, and the HRE case $p=2$    settled in the present paper; this latter case  is the most challenging, as can  be seen also from the proof in \cite{BK}. For other extensions of the \textsc{Baum-Katz} theorem to martingales (respectively, exchangeable  random variables), see   \cite{MYS} (resp.,\,\cite{St2}).  There is also a wide literature on limit theorems involving complete convergence even beyond the \textsc{Baum-Katz} type, but we shall not discuss this aspect in the present paper.

%%%%%%%%%%%%%%%%%%%%%%%
\section{The Proof of Proposition \ref{H-HRE-1}}
\label{sec4}
%%%%%%%%%%%%%%%%%%%%%%%
 
 The proof of Proposition \ref{H-HRE-1}  will involve several steps, which we have tried   to outline   clearly. 

\smallskip
 On the strength of assumption  $(\star)$, the sequence of functions $\, f_1, f_2, \cdots\,$ contains a (relabelled) subsequence, such that $\, f_1^2, f_2^2, \cdots\,$ converges weakly in $\Ll^1$ to $\, \boldsymbol{\eta}  \in \Ll^\infty$ as in \eqref{2.2}. Passing to a further subsequence, still denoted $\, f_1, f_2, \cdots\,,$ we may assume that   this sequence

 \noindent
{\bf .}   is {\it determining,} i.e., satisfies \eqref{2.7} for some limit random probability distribution   $\,   \mathbf{ H}  (\cdot\,, \omega)$ with first and second moments $f_\infty \in \Ll^2$ and $\, \boldsymbol{\eta}  \in \Ll^\infty$,  as in \eqref{2.8},

 \noindent
{\bf .} has squares $\, f_1^2, f_2^2, \cdots\,$ which are uniformly integrable as in \eqref{2.6} of subsection \ref{rem2.4}, and 

 \noindent
{\bf .}  {\it converges weakly in $\Ll^2$ to} $\,f_\infty \in \Ll^2\,,  $ i.e., 
 \begin{equation} 
\label{4.1}
\lim_{n \to \infty} \E \big( f_{ n}  \cdot \xi \big) = \E \big( f_\infty \cdot \xi \big)\,, \quad \forall ~~ \xi \in \Ll^2\,.   
\end{equation} 
These properties are shared by all further subsequences of 
 $\, f_1, f_2, \cdots\,$, as well as by their permutations. Without   sacrificing generality,     we take 
 \begin{equation} 
\label{4.1z}
  f_\infty \equiv 0\,,\quad  \big\|   \boldsymbol{\eta}  \big\|_\infty \le 1\,.   
\end{equation} 
  $\bullet\,$   We may also assume, in accordance with the   perturbation arguments of \textsc{Chatterji} (\cite{Ch}, pp.\,137-141; also         discussed    here in an Appendix, section \ref{sec7}), and   passing inductively to a subsequence  if necessary,  that the  $f_1, f_2, \cdots$ are {\it simple, $\Ll^2-$bounded     martingale differences}  \,(``strongly orthogonal" in the terminology of  
\cite{Ch3}), i.e., satisfy
\begin{equation} 
\label{4.3}
\E\big( f_{n+1} \, \big| \, {\cal F}_n \big) =0\,, ~~~\E\big( f_{n+1}^2 \, \big| \, {\cal F}_n \big) \le 1\,,  ~~\P-\text{a.e.,}\quad \text{with}\quad \F_n := \sigma \big(f_1, \cdots, f_n \big)
\end{equation} 
 for every $n \in \N$ and with $\, \F_0 := \{ \emptyset, \Omega\}$. Thus, the sequence  $\,X_N:= \sum_{n=1}^N f_n\,,~ N \in \N\,$ is then   a   square-integrable martingale  of the  resulting  filtration 
\begin{equation} 
\label{4.3a}
\mathbb{F}\,:= \big\{ \F_n \big\}_{n \in \N_0} \,,
\end{equation} 
with $\,\E \big( X_N^2\big) \le N$.   Since $f_\infty \equiv 0\,,$ it is enough  to establish the HRE property \eqref{2.4} with $\eps =1$ (because we can then replace the $f_1, f_2, \cdots$ by $f_1 / \eps, f_2/  \eps, \cdots)$;  i.e.,   show that  the sets 
\begin{equation} 
\label{4.7}
 B_N := \Big\{ \big| X_N \big| > N \Big\} = \bigg\{ \Big|  \sum_{n=1}^N f_n   \Big| > N \bigg\}\,, ~ N \in \N  ~  \quad \text{satisfy} \quad \sum_{N \in \N} \, \P \big( B_N    \big)   < \infty\,.
 \end{equation}
 $\bullet\,$ 
We start by introducing the stopping times
\begin{equation} 
\label{4.20}
\, \boldsymbol{\tau}_N := \min \Big\{ n=1, \cdots, N : \big| X_n \big| > N/3 \Big\}\,, \qquad N \in \N
\end{equation}
of the filtration $ \mathbb{F}  $ in \eqref{4.3a}   (with the understanding $\min \emptyset \equiv \infty$), so that the \textsc{\v Ceby\v sev}   and  \textsc{Doob}-maximal  (e.g.,\,\cite{Du}, pp.\,249-250) inequalities give
$$
\P \big(   \boldsymbol{\tau}_N < \infty \big)    
= \, \P \Big( \max_{1 \le n \le N} \big| X_n \big| > \frac{\,N\,}{3} \Big)
~~~~~~~~~~~~~~~~~~~~~~~~
$$  
\begin{equation} 
\label{4.21}
~~~~~~~~~~~~~~~~~~~~~~~~~~~~~~~~~ 
\le \frac{9}{\,N^2\,}\, \,\E\Big[ \, \max_{1  \le n \le N}   X^2_n  \, \Big] \le \frac{\,36\,}{\,N^2\,}\, \,\E \big(     X_N^2 \big) \le \frac{\,36\,}{\,N\,} \,.
\end{equation}
Summation over $N $ yields a   divergent series, so we need an improved estimate. In particular, we shall try to show that  $\big| X_{  \boldsymbol{\tau}_N}\big|$ is bounded by $2N/3$ with high probability; we shall do this  by estimating the possible size of the jump of $X_1, X_2, \cdots$ at time $\,\boldsymbol{\tau}_N\,$. 

We note that \eqref{4.21} still holds for all further subsequences of $\,f_1, f_2, \cdots\,,$ as well as for their permutations. The arguments that follow will also be designed   to share these hereditary features. More precisely, we shall find a further subsequence, still denoted $\,f_1, f_2, \cdots\,$, such that, for all permutations $\, f_{k_1},  f_{k_2},\cdots\,$  of this  (and of every further) subsequence, we have
\begin{equation} 
\label{4.6z}
\sum_{N \in \N} \, \P \,\bigg( \max_{ 1\le n \le N} \big| f_{k_n} \big|   > \frac{N}{3}\,  \bigg) < \infty\,.  
\end{equation}
 $\bullet\,$ To prove \eqref{4.6z}, we prepare the ground by applying the reasoning of \textsc{Erd\H{o}s} on p.\,287 in \cite{Erd}, developed there for the I.I.D. case,   but adapted  here  to our own situation of the more general sequences $\,f_1, f_2, \cdots\,$ as above. {\it We shall only use the $\Ll^2-$boundedness of such a sequence, the uniform integrability of $\,f_1^2, f_2^2, \cdots\,,$ and the assumptions of \eqref{4.1z}.}
 
 \smallskip
Let us start this preparation by observing that, for   $f  \in \Ll^2\,$ and with   $\, a_i := \P \big( |f | > 2^i\big)\,,\, \,i \in \N_0\,,$ we have
$$
\sum_{i \in \N_0} 2^{\,2i - 1}   a_i   \, \le  \, \sum_{i \in \N_0} 2^{\,2i } \big( a_i- a_{i+1}\big) \, \le  \, \E(f^2) \,  \le  \,  \sum_{i \in \N_0} 2^{\,2 (i+1) } \big( a_i- a_{i+1}\big) \, \le \, \sum_{i \in \N_0} 2^{\,2i +2}   a_i\,. 
$$
 These   inequalities show  that the second-moment condition $\, \E(f^2 )< \infty\,$  is equivalent to 
\begin{equation} 
\label{ERD2}
\sum_{i \in \N_0} 2^{\,2i  } \,  a_i < \infty\,.
\end{equation} 
Now, for the sequence   $\,f_1, f_2, \cdots\,$ with the properties stated above, we can assume   that 
\begin{equation} 
\label{ERD3}
    a^n_i := \P \big( |f_n | > 2^i \big)\,, \quad n \in \N
    \end{equation} 
     converges to some limit $\,a_i\,$, as $\, n \to \infty\,,$ for each $\, i \in \N_0\,;$ and by $\Ll^2-$boundedness, we have   the property \eqref{ERD2} for these limits. Whereas, by passing once more to a (relabelled) subsequence   again, we obtain 
     for each $\, i \in \N\,$ the bound
\begin{equation} 
\label{ERD6}
 a^n_i\,<\,a_i + 2^{\,- 3i}\,, \qquad \forall ~~n \ge i\,.
    \end{equation}    
      {\it   Pretending for a moment that  \eqref{ERD6}    holds for all $\, (n,i) \in \N \times \N_0\,,$} we  establish  the estimate
       $$
   \sum_{   N \,= \,2^{\, i}}^{ 2^{\, i+1}-1}    \,\P \,\Big( \max_{1 \le n \le N} \big| f_{ n} \big|   >   N\Big) \, \le\, \sum_{   N \,= 2^i}^{ 2^{\, i+1} -1} \,\sum_{n=1}^N    \,\P \big( |f_n | > 2^i \big)
   $$
   $$
 ~~~~~~  \, \le \,   2^i  \, \sum_{n=1}^{2^{\, i +1 } }  a^n_i\,<\, 2^{\,2i +1 } \, \big(  a_i + 2^{\,-3 i } \big)\,, \quad \forall~ i \in \N_0\,,
  $$ 
  and   note that the right-hand-side is summable over $i \in \N$ on account of \eqref{ERD2}. 
   
   Of course, for given, fixed $i \in \N_0\,$, we cannot guarantee the validity of \eqref{ERD6} for all $\,n \in \N\,;$ but we {\it can} obtain a (relabelled) subsequence of $\,f_1, f_2, \cdots\,$ which satisfies, for some $ i_0 \in \N_0\,,$ the bound
\begin{equation} 
\label{ERD5}
    a^n_i  = \P \big( |f_n | > 2^i \big)\,<\,2^{\,- 2i}\,, \qquad \forall ~~n \in \N
    \end{equation}    
 for all integers  $\,i \ge i_0\,;$ for otherwise the presumed uniform integrability of the    $\,f_1^2, f_2^2, \cdots\,$ would fail.  Combining the two estimates \eqref{ERD6}, \eqref{ERD5} we deduce that  a further subsequence, again denoted $f_1, f_2, \cdots\,,$   can be selected, along which   \eqref{ERD5} holds   for {\it  every} $\,i \in \N_0\,.$  
 
Now, for each $i \in \N_0$ and $N \in \N$   with   $\, 2^i \le N < 2^{\,i+1}\,,  $ there are at most $i$ terms in the string $\,  1,  2, \cdots,  N\, $ for which \eqref{ERD6} can fail. For these ``offending" terms we can apply \eqref{ERD5} rather than \eqref{ERD6}, and obtain (without using the above italicized ``pretense") the upper bound 
 \begin{equation} 
\label{ERD7}
 \sum_{   N \,= \,2^{\, i}}^{ 2^{\, i+1}-1}    \,\P \,\Big( \max_{1 \le n \le N} \big| f_{ n} \big|   >   N \Big) \, \le\,   2^{\,2 i + 1} \Big( a_i + 2^{\,- 3 i} \Big) +2^i \cdot i \, 2^{\,- 2i}\,, \quad i \ge i_0\,,
    \end{equation}
which is again the general term of a   convergent series. 

\smallskip
\noindent 
 $\bullet\,$ Returning to \eqref{4.6z}, we choose any permutation $\,  ( k_n  )_{n \in N}\,$ of a subsequence of the above-constructed $\, f_1, f_2, \cdots\,.$ In other words, we fix any sequence $\,  ( k_n  )_{n \in N}\,$ of distinct natural numbers; and have to show that the ``perturbed sequence" $\, f_{k_1},  f_{k_2},\cdots\,$ has the HRE property, thus finishing the proof of Proposition \ref{H-HRE-1}. To wit, with the notation $\, \widetilde{X}_N := f_{k_1} + \cdots + f_{k_N}\,,$
 $$
 \widetilde{A}_N := \bigg\{\, \max_{1\le n \le N} \big| f_{k_n} \big|   \le   \frac{N}{3}\, \bigg\} \,, \quad \widetilde{B}_N := \Big\{ \big| \widetilde{X}_N \big| > N \Big\} = \bigg\{ \Big|  \sum_{n=1}^N f_{k_n}   \Big| > N \bigg\}\,, \quad N \in \N
 $$
we have to show $ \, \sum_{N \in \N} \, \P ( \widetilde{B}_N ) < \infty\,,$ as in \eqref{4.7}.  But, in the light of \eqref{4.6z}, it suffices to   prove  
\begin{equation} 
\label{4.10z}
\sum_{N \in \N} \, \P \big( \widetilde{B}_N \cap \widetilde{A}_N  \big)   < \infty\,.
\end{equation}
To estimate this probability, and denoting by $\,\widetilde{\boldsymbol{\tau}}_N := \min \big\{ n=1, \cdots, N : \big| \widetilde{X}_n \big| > N/3 \big\}\, $ the analogue of the stopping time in \eqref{4.24}, we observe that $\om \in \widetilde{B}_N$ implies $\widetilde{\boldsymbol{\tau}}_N (\om) < \infty$ and $\big| \widetilde{X}_{ \widetilde{\boldsymbol{\tau}}_N (\om)   -1} (\om) \big| \le N/3\,.$ The jump $\, \widetilde{X}_{ \widetilde{\boldsymbol{\tau}}_N (\om)   } (\om) - \widetilde{X}_{ \widetilde{\boldsymbol{\tau}}_N (\om)   -1} (\om)\,$ equals $f_{ \widetilde{\boldsymbol{\tau}}_N (\om)   } (\om)$ and thus, for $\,\om \in A_N$, we have $\, \big| \widetilde{X}_{ \widetilde{\boldsymbol{\tau}}_N (\om)   } (\om) \big| \le 2N/3\,.$ 

Finally, we  need to control the probability that $\, \big| \widetilde{X}_N - \widetilde{X}_{ \widetilde{\boldsymbol{\tau}}_N} \big|$ can become bigger than $ N / 3$. For this purpose, we introduce the    $\Ll^2-$bounded  martingale
\begin{equation} 
\label{4.24}
\widetilde{Y}_n \,:=\,  \widetilde{X}_n -\widetilde{X}_{n \wedge   \widetilde{\boldsymbol{\tau}}_N } \,=\, \Big( \widetilde{X}_n -\widetilde{X}_{   \, \widetilde{\boldsymbol{\tau}}_N } \Big) \cdot \Ind_{ \{    \widetilde{\boldsymbol{\tau}}_N  < n  \} } \,, \qquad n=1, \cdots, N 
\end{equation}
with $\, \E \,\big[ \big( \widetilde{Y}_n-\widetilde{Y}_{n-1} \big)^2 \, \big|\, {\cal F}_{ \, \widetilde{\boldsymbol{\tau}}_N}\big] \le 1\,, $ and thus also $  \P \,  \big( \big| \widetilde{Y}_N \big| > N/3\, \big|\, {\cal F}_{\, \widetilde{\boldsymbol{\tau}}_N}\big) \le 9/N\,$ as in \eqref{4.21}. 
This gives the desired control on the difference $\, \widetilde{Y}_N = \widetilde{X}_N - \widetilde{X}_{\widetilde{\boldsymbol{\tau}}_N}\,$, conditionally on ${\cal F}_{ \, \widetilde{\boldsymbol{\tau}}_N}$. 

Indeed, on the event $\big\{ \widetilde{\boldsymbol{\tau}}_N = n \big\}\,,$ we may estimate now 
$$
\P \big( \widetilde{B}_N \cap \widetilde{A}_N , \,\widetilde{\boldsymbol{\tau}}_N = n \big) \, \le \, \P \big( \widetilde{Y}_N > N / 3, \, \widetilde{\boldsymbol{\tau}}_N = n \big) \, \le \, \frac{9}{N} \cdot \P \big( \widetilde{\boldsymbol{\tau}}_N = n \big)\,;
$$
whereas,  summing up over $n=1, \cdots, N$  and recalling  \eqref{4.21}, we obtain the bound 
$$
\P \big( \widetilde{B}_N \cap \widetilde{A}_N   \big)\,\le\, \frac{9}{N} \cdot \sum_{n=1}^N \, \P \big( \widetilde{\boldsymbol{\tau}}_N = n \big) \,\le\, \frac{9}{N} \cdot \frac{36}{N} \,,
$$
which yields the desired convergent series as in \eqref{4.10z}. \qed

%%%%%%%%%%%%%%%%%%%%%%%
\section{Strong Exchangeability at Infinity; Independence}
\label{sec6}
%%%%%%%%%%%%%%%%%%%%%%%

We prepare now the ground for the proof of Theorem \ref{H-HRE}\,(i),  following the trail blazed in the seminal works  of \textsc{Aldous} \cite{A}\,--\,\cite{AE}  and,  significantly, \textsc{Berkes-P\'eter} \cite{BerPet}.  Our  path  here will be straighter, as we can work   in an $\Ll^2-$setting and   thus   measure distances between probability measures using the quadratic \textsc{Wasserstein},  rather than the more delicate \textsc{Prokhorov}, metric.  This allows simpler arguments which we have endeavored, nevertheless, to spell out in     detail.

 \smallskip
We   assume   as in \cite{BerPet}     that the underlying probability space $(\Omega, \F, \P)$ is  separable,  and rich enough to accommodate all the objects we   need to place on it: in particular,   a sequence $f_1, f_2, \cdots$ of functions   as in  Theorem \ref{H-HRE}\,(i),    for which   we may assume        $\,f_1^2 \,\mathbf{ 1}_{A_{1}}, f_2^2 \, \mathbf{ 1}_{A_{2}}, \cdots$ to be  uniformly integrable    in   (cf.\,\eqref{2.6}).  As in  section \ref{sec4}, we       take   $\,f_1  \,\mathbf{ 1}_{A_{1}}, f_2  \, \mathbf{ 1}_{A_{2}}, \cdots$ to be      simple,   martingale-differences  (cf.\,\eqref{4.3});   and observe from    \eqref{4.1}  that, after  passing to a (re-labelled) subsequence, we may assume also in the manner of \eqref{4.1},    for some $f_\infty \in \Ll^2\,$, 
\begin{equation} 
\label{6.1}
\lim_{n \to \infty} \E \big( f_{ n}\, \mathbf{ 1}_{A_{n}}  \cdot \xi \big) = \E \big( f_{\infty}    \cdot \xi \big)\,, \quad \forall ~~ \xi \in \Ll^2 .
\end{equation} 
To simplify typography we shall  take    $f_\infty \equiv 0\,$  
 and  denote, in the rest of this section  and in section \ref{sec5}, this simple,   martingale-difference sequence $\,f_1  \,\mathbf{ 1}_{A_{1}}, f_2  \, \mathbf{ 1}_{A_{2}}, \cdots$          simply by $\,f_1, f_2, \cdots \,$.

 \smallskip
\noindent
$\bullet\,\,$ 
We recall  now  the filtration $ \mathbb{F}  =  \{ \F_n  \}_{n \in \N_0} $ from \eqref{4.3a}; consider   atoms  $\, A^{(n)}_j,$ $ j=1,\cdots, J_n\,$    of positive $\P-$measure  which generate  the $\sigma-$algebra $\F_n\,$ for each   $n \in \N\,$; and denote the conditional distribution of $f_n\,,$ given a generic  one of these atoms $A$ with $\P(A)>0$, by 
\begin{equation} 
\label{6.2}
\mu_{(n)}^A (\cdot) \,:=\, \P \big( f_n \in \,\cdot \,\,\big| \, A \big)\,.
\end{equation} 
 This    is an element of    ${\cal P}_2 (\R), $ the space of probability measures on the \textsc{Borel} sets of the real line   with  finite second moment. We  endow this  ${\cal P}_2 (\R)$ with the quadratic \textsc{Wasserstein} distance below, which renders it  a  Polish (i.e.,   complete, separable, metric;   cf.\,\cite{Vill})  space:
$$
{\cal W}_2 \big( \mu, \nu) \,:= \, \left( \inf_{X \sim \mu\,,\, Y \sim \nu} \E \Big( X-Y\Big )^2 \, \right)^{1/2} .
$$

 Now,    the assumed uniform integrability    of the sequence $\,   f_1^2, f_2^2, \cdots    \,$     implies that, for each fixed atom $A\,$ in $\big( A^{(n)}_j \,;\,j=1,\cdots, J_n \,, \, n \in \N\, \big) ,$  the sequence of probability measures
\begin{equation} 
\label{6.3}
{\cal M}_A \,:=\, \big\{ \mu_{(n)}^A    \big\}_{n \in \N}
\end{equation} 
in \eqref{6.2} is {\it tight:}   given any $\varepsilon >0\,,$ there exists a compact     $K_\varepsilon \subset \R$     with   $\, \int_{\,\R \setminus K_\varepsilon} x^2 \, \nu ( \mathrm{d} x) < \varepsilon\,,$ $  \, \forall ~\, \nu \in {\cal M}_A \,$. Thus $\,{\cal M}_A \, $ is  a {\it relatively compact} subset of    ${\cal P}_2 (\R)$:  given any   atom $A$, we can find a probability measure $\,\mu^A \in {\cal P}_2 (\R)\,  $  and   a    subsequence $f^{A}_{k_1}, f^{A}_{k_2}, \cdots$ of $f_1, f_2, \cdots,$ with
 \begin{equation} 
\label{6.4}
\lim_{n \to \infty} {\cal W}_2 \big( \mu^A, \, \mu_{(k_n)}^A\big)\,=\,0 
\end{equation}
 (Theorem 7.12 in \cite{Vill});  and eventually, by diagonalization, also an "omnibus" subsequence $\,f_{k_1}, f_{k_2}, \cdots\,$ of $\,f_1, f_2, \cdots \,$ with     \eqref{6.4}   valid for {\it each} such atom $A$.

 This $\,f_{k_1}, f_{k_2}, \cdots \,$ is a martingale-difference sequence  of the "thinned" sub-filtration $\big\{ \F_{k_n}\big\}_{n\in \N}$ of   
$\big\{ \F_n\big\}_{n\in \N}$   in \eqref{4.3a}.   To ease    (the already   heavy) notation somewhat, we   denote these thinned subsequences as  $\big\{ f_{ n} \big\}_{n \in \N}\,$, $\big\{ \mu_{ (n)}^A \big\}_{n \in \N}$  and $\big\{ \F_n\big\}_{n\in \N}\,$, respectively.   
 
 \smallskip 
 \noindent 
 $\bullet \,$ 
We follow then again the trail from \textsc{R\'enyi}\,\cite{Ren} (cf.\,\cite{AE},\,\cite{BerRos})  and obtain the existence of a measurable mapping $\,   \boldsymbol{\mu}  : \Omega \to {\cal P}_2 (\R)\,$ which     aggregates    the limiting probability measures in \eqref{6.4}, in the sense that,  for {\it each} atom $A$  as above, we have  the property 
\begin{equation} 
\label{6.5}~ 
\mu^A (\cdot) = \int_A   \boldsymbol{\mu}  (\cdot \, , \omega) \, \P (\mathrm{d} \omega)\,, ~~ \text{i.e.,} ~~
\int_\R \varphi (x) \, \mu^A (\mathrm{d} x) \,= \int_A \left( \,\int_\R \varphi (x)\,    \boldsymbol{\mu} ( \mathrm{d} x, \omega) \right) \P (\mathrm{d} \omega)
\end{equation}
for every bounded and continuous $\varphi : \R \to \R\,$. 
 This aggregating probability measure  $   \boldsymbol{\mu}  (\cdot\,, \omega)$ has   $\,x \mapsto   \mathbf{ H}  (x, \omega)$ of \eqref{2.7} as its distribution function.  We define also, for each fixed $n \in \N$, an $\F_n-$measurable mapping $\,    \boldsymbol{\mu}_{(n)} : \Omega \to {\cal P}_2 (\R)\,$ via
 \begin{equation} 
\label{6.5a}
   \boldsymbol{\mu}_{(n)} (\omega, \cdot) \,:=\, \mu_{(n)}^A (\cdot) \,, \qquad \omega \in \Omega
 \end{equation} 
as in \eqref{6.2},   with $A$ the unique atom        among   $\, A^{(n)}_j,~ j=1,\cdots, J_n\,$ for which $ \,  \omega \in A\,.$  In this manner we obtain  the aggregated (and relabelled) version on \eqref{6.4}, namely, 
\begin{equation} 
\label{6.6}
\lim_{n \to \infty} {\cal W}_2 \big(   \boldsymbol{\mu} \,,    \boldsymbol{\mu}_{( n)} \big)\,=\,0\,, ~\P-\text{a.e.;~~~ as well as} ~~~~\lim_{n \to \infty} \, \E \,\big[ \,  {\cal W}_2 \big(    \boldsymbol{\mu} \,,    \boldsymbol{\mu}_{( n)} \big)\big] \,=\,0\,.
\end{equation}
Furthermore, appealing to \textsc{Egorov}'s theorem, we obtain sets $\, E_1 \subseteq  E_2 \subseteq  \cdots\,$ in $\F$ with $\, \P (E_m) \ge 1 - \varepsilon_m\,$ and $\lim_{m \to \infty} \downarrow \varepsilon_m=0\,$, such that the $\P-$a.e.\,convergence in \eqref{6.6} is uniform on each $E_m\,$. We may   assume  additionally   that, for each of these sets $E_m\,$,  the restriction $   \boldsymbol{\mu} \big|_{E_m}\,$ of  the aggregator $    \boldsymbol{\mu}  \,$ in \eqref{6.5}  is supported on a relatively compact subset of $\,{\cal P}_2 (\R).$  

\smallskip
   We  adapt  now  the    notion of   strong exchangeability at infinity  from    \cite{BerPet}       to our $\Ll^2-$\,setting.

\begin{definition} 
 {\bf Strong $\Ll^2-$Exchangeability at Infinity.} 
\label{def_sei}
Fix a sequence of    functions $\, \big( g_n \big)_{n\in \N} \subset \Ll^2$, and a sequence $\,   \boldsymbol{\varepsilon}  :=\big\{ \varepsilon_k \big\}_{k \in \N} \subset (0,1)\,$ with $\,\lim_{k \to \infty} \downarrow \varepsilon_k=0\,.$ 

We call  $\, \big( g_n \big)_{n\in \N}\,$     {\it strongly $\Ll^2-$exchangeable at infinity  with speed $  \boldsymbol{\varepsilon} $} if, {\rm for  each} $k \in \N$, there exists    a partition $\, \big\{ A^{(k)}_0, A^{(k)}_1, \cdots, A^{(k)}_{J_k}   \big\} \subset \F\,$ of $\,\Omega$ by disjoint sets    of positive $\P-$measure  with  the following  properties:

\noindent
{\bf (i)}   $\,A^{(k+1)}_0 \subseteq A^{(k )}_0\,,$ and $\, \big( A^{(k+1)}_j,~ j=0,1,\cdots, J_{k+1} \big)\,$    refines $\, \big( A^{(k)}_j,~ j=0,1,\cdots, J_k \big)\,;$

\noindent
{\bf (ii)}  
$\,
\P \big( A^{(k )}_0 \big) \le \varepsilon_k\,,~~ \sup_{n \in \N} \,\E \big( g_n^2 \cdot \Ind_{A^{(k )}_0} \big) \le \varepsilon_k\,; 
 $

  \noindent
{\bf (iii)} For each  $A$ among the $\,  A^{(k)}_1, \cdots, A^{(k)}_{J_k}  \, ,$ there exist   independent, square-integrable functions 
$\,
h^A_{k+1}\,, h^A_{k+2}\,,$ $ \cdots 
\,,$
 with  common distribution $\mu^A$ as in \eqref{6.4} and the sentence preceding it, and with the property  
$\, 
\E^{\,\P^A } \big[\, \big( g_n -h^A_{n}  \big)^2\,\big] \le \,\varepsilon_k\,, ~\,  \forall ~ n \ge k+1\,  
 $ under the conditional probability measure 
\begin{equation} 
\label{6.7}
\P^A \big(\cdot \big) \,:=\,  \P  \big(\cdot \,\cap \,A \, \big)\, / \,\, \P \big(A \big)  \,.
\end{equation} 
\end{definition}

Here is    an analogue of Theorem 1 in  \cite{BerPet},    tailored to our $\Ll^2-\,$setting. It approximates (subsequences of) martingale-difference sequences by exchangeable ones, in a "progressively improving" manner.

\begin{theorem} 
\label{thm1b}
Fix a sequence of   
functions $\, f_1, f_2, \cdots\,$  as in   this section;  and a  decreasing sequence $\, \boldsymbol{\varepsilon}:=\big\{ \varepsilon_k \big\}_{k \in \N} \subset (0,1)\,$ with $\,\lim_{k \to \infty} \downarrow \varepsilon_k=0\,.$ 

There exist then a subsequence $\, f_{\ell_1}, f_{\ell_2}, \cdots\,$ and, {\rm for each} $\,k \in \N\,,$    {\rm exchangeable} square-integrable functions       $\,\widehat{h}^{(k)}_{k+1},\, \widehat{h}^{(k)}_{k+2}, \cdots\,$, so that   we have 
\begin{equation} 
\label{6.7FF}
 \E \big[ \,\big( f_{\ell_n} -  \widehat{h}^{(k)}_n \big)^2  \, \big|\, \F_k \, \big] \le \varepsilon_k\,, \quad \forall ~ n \ge k+1  \quad ~\text{with} \quad \F_k \,:=\, \sigma \big( f_{\ell_1}, f_{\ell_2}, \cdots,  f_{\ell_k}\big)  \,.
 \end{equation} 
\end{theorem}

In a similar spirit, we formulate an analogue of Theorem 2 in \textsc{Berkes-P\'eter}\cite{BerPet}, dealing with strong exchangeability and  tailored once  again  to our $\Ll^2-$setting. Just as in   \cite{BerPet} (whose   Theorem 1 follows from   Theorem 2 there),  Theorem \ref{thm1b} right above is a direct consequence of our next result, Theorem \ref{thm2b}; this  is    proved  in subsection \ref{sec4.2} below.

\begin{theorem} 
\label{thm2b}
Fix a sequence   $\, f_1, f_2, \cdots\,$  as in   this section;  and a   decreasing sequence $\,   \boldsymbol{\varepsilon}  :=\big\{ \varepsilon_k \big\}_{k \in \N} \subset (0,1)\,$ with $\,\lim_{k \to \infty} \downarrow \varepsilon_k=0\,.$

There exist then   a subsequence   $\, f_{\ell_1}, f_{\ell_2}, \cdots\,$ of $\,f_1, f_2, \cdots\,,$ and    a sequence $\, \big( g_n \big)_{n \in \N} \subseteq \Ll^2\,$ strongly $\,\Ll^2-$exchangeable at infinity    with speed   $ \boldsymbol{\varepsilon}  ,$       so that  the analogue of \eqref{6.7FF}, namely,   $\,\, \E \big[ \big( f_{\ell_n} - g_n \big)^2\, \big|\, \F_k \big] \le \varepsilon_k\,,$ holds for   $\,(n,k) \in \N^{\,2} $   with $\, n \ge k+1  \,.$ 
\end{theorem}

We deduce in the next subsection some important consequences of this result. We start by   casting  it in an equivalent but  more detailed and operational   form, recalling  Definition \ref{def_sei}.

\begin{corollary} 
\label{cor2b}
In the setting of Theorem \ref{thm2b}, there exist a   subsequence   $\, f_{\ell_1}, f_{\ell_2}, \cdots\,$ of the given sequence $\,f_1, f_2, \cdots\,,$ and a double array $\,   \big( A^{(k)}_j,~ j=  0,1,\cdots, J_k \big)_{k \in \N}\,$ of row-wise disjoint sets  with  $\,\P \big( A^{(k)}_j \big) >0\,,$   such that, {\rm for each given} $\,k \in \N\,,$ we have:

 \smallskip
\noindent
{\bf (i)}   $\,A^{(k+1)}_0:= \Omega \setminus \bigcup_{j=1}^{J_{k+1}} A^{(k+1)}_j  \subseteq \Omega \setminus \bigcup_{j=1}^{J_k} A^{(k)}_j =: A^{(k )}_0\,,$ and  the partition 
$\, \big( A^{(k+1)}_j,~ j=0,1,\cdots, J_{k+1} \big)\,$ is a refinement of the preceding partition $\, \big( A^{(k)}_j,~ j=0,1,\cdots, J_k \big)\,;$

\noindent
{\bf (ii)}  
$\,
\P \big( A^{(k )}_0 \big) \le \varepsilon_k\,;~~~ \sup_{n \in \N} \,\E \Big( f_{\ell_n}^2 \cdot \Ind_{A^{(k )}_0} \Big) \le \varepsilon_k\,; 
 $
 
  \noindent
 {\bf (iii)} the $\sigma-$algebra $\F_k$ of \eqref{6.7FF} is included in $\, {\cal G}_k = \sigma \big( A^{(k)}_j\,, ~ j=0, 1, \cdots , J_k \big)\,;$ and

 \noindent
{\bf (iv)} for each set $A$ among the $\,  A^{(k)}_1, \cdots, A^{(k)}_{J_k}  \, ,$ there exist  square-integrable functions 
$\,
h^A_{k+1}\,, h^A_{k+2}\,, \cdots 
\,,$
 independent  and with      common distribution  $\mu^A$  as in  \eqref{6.4}   under the   probability measure $\,\P^A   $ of  \eqref{6.7}, which satisfy 
$\,~
\E^{\,\P^A} \big( f_{\ell_n} - h^A_{n} \big)^2 \le \,\varepsilon_k\,, ~  \forall \,\, n \ge k+1\,; 
 $ i.e.,
 \begin{equation} 
\label{6.7a}
\E \,\Big[\, \Big( f_{\ell_n} - h^{A^{(k)}_j}_n  \Big)^2  \Ind_{A^{(k)}_j} \, \Big] \le \,\varepsilon_k  \, \P\big(  A^{(k)}_j   \big)\,, \quad  \forall ~ j=1, \cdots, J_k\,,~~n = k+1, \,k+2, \cdots\,.
\end{equation} 
 \end{corollary}

%%%%%%%%%%%%%%%%%%%%%%%
\subsection{Approximation by an Omnibus Sequence}
\label{sec4.1}
%%%%%%%%%%%%%%%%%%%%%%%

Corollary \ref{cor2b} casts Theorem \ref{thm2b} in terms of a "progressively improving"  (i.e.,  with diminishing error $\,   \varepsilon_k \downarrow0$) $\,\Ll^2-$approximation of an appropriate subsequence   $\, f_{\ell_1}, f_{\ell_2}, \cdots\,$ of $f_1, f_2, \cdots$, using double  arrays $\big( \,  \widehat{h}^{(k)}_{k+1},\, \widehat{h}^{(k)}_{k+2}, \cdots\,  \big)_{k \in \N } \,$ of row-wise independent functions in $\Ll^2$ with common distribution.   

This setup is almost exactly that of       \cite{BerPet}; it imposes only  `tightness'  on  the $f_1, f_2, \cdots,$ and     uses \textsc{Prokhorov}  distances. But   it comes at a price: {\it at each level} $k \in \N\,,$ it has to  start a {\it new}  register   $k+1, \,k+2, \cdots \,,$ and discard an exceptional set $\, A^{(k )}_0\, $ of small $\P-$measure.   

 \smallskip
 
For our purposes,  we shall need only $\Ll^2-${\it bounded},   as apposed to $\Ll^2-$small-and-diminishing, approximations of the terms in the subsequence   $\, f_{\ell_1}, f_{\ell_2}, \cdots\,.$  The fact that we have already proved Proposition \ref{H-HRE-1} will afford  us this small luxury;  and  will enable us to put together,    instead of a double array $\big( \,  \widehat{h}^{(k)}_{k+1},\, \widehat{h}^{(k)}_{k+2}, \cdots\,  \big)_{k \in \N } \,$ as in Corollary \ref{cor2b},  {\it a single, omnibus  sequence} $\, h_1, h_2, \cdots\,$ of independent, centered  functions,

\noindent
 {\bf (a)} whose   conditional distributions, given each set in  a partition of $\,\Omega$, are the same;   and \\ {\bf (b)} which approximates, in a good $\Ll^2-$sense, a   sequence    $\,f_1^*, f_2^*, \cdots\,$  suitably close to the subsequence $\, f_{\ell_1}, f_{\ell_2}, \cdots\,.$\,

 In this manner    we will not have   to restart a new exchangeable   sequence $\, \widehat{h}^{(k)}_{k+1},\, \widehat{h}^{(k)}_{k+2}, \cdots\,   $   at each level  $\,k \in \N$ of approximation.

  \smallskip
  \noindent
  $\bullet~$
 We put together  now this omnibus  sequence $\, h_1, h_2, \cdots\,$. We construct  the subsequence   $\, f_{\ell_1}, f_{\ell_2}, \cdots\,$ and  the partitions   $\, {\cal A}^{(k)}= \big( A^{(k)}_j,~ j=  0,1,\cdots, J_k \big)\,,$ ${k \in \N}\,$ by induction:
  \\ {\bf (i)}\, Suppose that   indices $\, \ell_1, \cdots, \ell_{k-1}\,$ and   atoms $ \, {\cal A}^{(\kappa)}= \big( A^{(\kappa)}_j,~ j=  0,1,\cdots, J_{\kappa} \big)\,, ~\kappa =1, \cdots, k-1\,,$ have been  selected.

 \noindent
  {\bf (ii)} \,Now define $\,\ell_k\,$ and $  {\cal A}^{(k)}\,$ by splitting the "exceptional atom" $\,A^{(k-1)}_0\,$ into disjoint sets  $\big( A^{(k)}_j,~ j=  0,1,\cdots, M_{k } \big)$  such that, for each given $\,j=1, \cdots, M_k\,$   and with $\,A \equiv   A^{(k)}_j\,,$ there are independent    $\,h^{A}_{\ell_k}, h^{A}_{\ell_k+1}\,, \cdots$ supported on $\,   A\,,$ with the same distribution,    and satisfying
 \begin{equation} 
\label{6.7b}
\big\| h^{A}_n - f_n \big\|_{\Ll^2 (\P^{A})} <1 \,, \quad \forall~ n \ge \ell_k\,.
\end{equation} 
{\bf (iii)}\, As for $A^{(k)}_0,$ we require $\, \P \big( A^{(k)}_0 \big) < 2^{-(k+1)}\,,$   $\,
 \big\|   f_n \cdot \Ind_{A^{(k)}_0}  \big\|_{\Ll^2 (\P )} <\,2^{-(k+2)} \,, ~~ \forall~ n \ge \ell_k\,;
  $
whereas, on the strength of this relation for the preceding step $k-1   $, we may   assume also 
 \begin{equation} 
\label{6.7c}
 \bigg\| \, \sum_{j=1}^{M_k} \,h^{A^{(k)}_j}_n \, \bigg\|_{\Ll^2 (\P )} < \, 2^{-(k+1)} \,.
\end{equation} 
{\bf (iv)} \,Regarding the remaining   atoms $\big( A^{(k-1)}_j,~ j=  1,\cdots, J_{k-1} \big)\,$ of the partition $\, {\cal A}^{(k-1)}\,,$ we {\it do not split them any further;} rather,   we keep them in the partition $\, {\cal A}^{(k)},$ {\it but relabelled,} namely as $\big( A^{(k)}_j,~ j=    M_{k }+1, \cdots, J_k \big)\,,$ so that $\, M_k + J_{k-1} = J_k\,.$ 
 
 \smallskip
Continuing in an obvious manner, we obtain a subsequence $\, f_{\ell_1}, f_{\ell_2}, \cdots\,$ as well as a countable partition ${\cal B}$ consisting of those atoms which appear  as  $\, A^{(k)}_j \in {\cal A}^{(k)}\,$ for some $(j,k)  \in \N^2\,$   (thus  also as   elements of $\,{\cal A}^{(m)}\,$ for $m>k$). Fixing an atom $B$ of this countable partition $\,   {\cal B}\,,$ we choose a   sequence $\, h^B_1, h^B_2, \cdots$ of independent and equi-distributed functions, supported by   $B$ and satisfying, with  $\, \boldsymbol{\kappa} (B)\,$   the smallest integer $k$ for which $B \in {\cal A}^{(k)}\,,$  the following bound, which  provides the desired estimate for integers $\,n \ge   \boldsymbol{\kappa}  (B)$:
 \begin{equation} 
\label{6.7d}
\big\| h^B_n - f_n \big\|_{\Ll^2 (\P^B)}  < 1 \,, \quad \forall~ n \ge   \boldsymbol{\kappa} (B) \,.
\end{equation}

But what about integers $n <   \boldsymbol{\kappa}  (B)$? To take care of these, we   modify  the sequence $\, f_{\ell_1}, f_{\ell_2}, \cdots\,,$ and obtain a new approximating sequence $\,f_1^*, f_2^*, \cdots\,$  as follows: For $B \in {\cal B}$ we define $\, f^*_n \, \Ind_B = h^B_n$ for $n <   \boldsymbol{\kappa}  (B)\,,$ so that    $\big\| h^B_n - f_n^* \big\|_{\Ll^2 (\P^B)} = 0<1 $ holds trivially.  Whereas, using \eqref{6.7b}--\eqref{6.7d},   the estimate $\, \E  [\, ( f_{\ell_k}- f_{k}^*   )^2 \, ] < 2^{-k}\,  $  holds.

\smallskip

In this manner, the sequence  $\, f_{\ell_1}, f_{\ell_2}, \cdots\,$ satisfies \eqref{2.4} with $f_\infty = 0$ (together with all its subsequences and permutations) if, and only if,    the so-constructed approximating sequence  $\, f_{1}^*, f_{2}^*, \cdots\,$ does.   We reprise   all the above by formulating  the central result of the present section,  an    omnibus version  of Theorem \ref{thm2b} and of Corollary \ref{cor2b}. 

\begin{theorem} {\rm An Omnibus $\Ll^2-$Approximation.} 
\label{prop2b}
For  a sequence of   
functions $\, f_1, f_2, \cdots\,$  as in this section, there exist  \\  {\bf .}  a subsequence   $\, f_{\ell_1}, f_{\ell_2}, \cdots\,$ of $\, f_1, f_2, \cdots\,;$  
\\  
{\bf .}  an approximating  sequence $\, \big( f_{ n}^* \big)_{n \in \N} \subset \Ll^2\,$    
with $\, \E \big[ \big( f_{\ell_n}- f_{ n}^*  \big)^2 \, \big] \le 2^{-n}\,$ for each $\,n \in \N\,;$ 
\\
{\bf .} a countable partition  $\, {\cal B}= \{ B_1, B_2, \cdots \}\,$ of  $\,\Omega$ by sets in $\F$ with positive $\P-$measure; and 
\\
 {\bf .}  square-integrable functions $\, h_1, h_2, \cdots \,$ which   are independent    and identically   distributed under the conditional probability measure $\,\P^{B}$ as in \eqref{6.7},  {\rm for each set} $\, B \in {\cal B}\,$, and satisfy  
\begin{equation} 
\label{6.8}
\E^{\,\P^{B}} \big[ \, \big( f^*_n - h_n \big)^2\,    \big]  \le 1\,,~~~ \text{that is\,}, ~~~ ~\E \Big[ \, \big( f^*_n - h_n \big)^2 \cdot \Ind_{B} \,   \, \Big] \le \,\P \left(    B   \right)\,, ~~~\forall ~n \in \N\,.
\end{equation} 
\end{theorem}

The  omnibus  sequence of square-integrable functions $\, h_1, h_2, \cdots \,$ in this Theorem  has properties  particularly well-suited to our context, as the following result demonstrates    (cf.\,Theorem 1 in \cite{St2}). 

\begin{proposition} 
\label{prop3}
Suppose the  functions  $\, f_1, f_2, \cdots \,$    are square-integrable, 
and that $\, {\cal B} = \{ B_1, B_2, \cdots \} $  is a partition of    $\,\Omega$ by sets in $\F$ of positive measure      such  that, conditioned on each $B_m\,,$ the  $\, f_1, f_2, \cdots \,$ are  independent and equistributed, with zero mean. 

  Then the $\, f_1, f_2, \cdots \,$  satisfy the HRE property \eqref{2.4} with $f_\infty \equiv 0\,$. 
  \end{proposition}

\noindent
{\it Proof:}  As already observed, it suffices to establish \eqref{2.4b} for $ \,\eps = 1\,$.  The uniform version of  the representation  \eqref{1.6}     in Proposition \ref{prop8.0},   provides   a universal constant $C>0$   such that  
$$  
\sum_{N \in \N} \P^{\,B_m}   \bigg( \bigg|  \sum_{n=1}^N f_n   \bigg| > N   \bigg) \le C \, \big( \sigma_m^2 \vee 1 \big) \,,  \quad \forall ~~ m \in \N
$$  
holds with $ \sigma_m^2 := \E^{\,\P^{B_m}} \big(f_1^2\big)$,  for $\, f_1, f_2, \cdots \,$ and  all its subsequences and permutations,  because of independence and common distribution under the probability measure $\P^{\,B_m}$. Multiplying by $\,\P(B_m),$ then summing up over $m \in \N,$ we obtain from the law of total probability  
$$
 \sum_{N \in \N} \P   \left( \Big|  \sum_{n=1}^N f_n   \Big| > N   \right) \le\,  C \sum_{m \in \N }  \P (B_m) \, \big(  \sigma_m^2  \vee 1 \big) \, \le \,  C \, \Big( 1 +\E^{\,\P} \big( f_1^2   \big) \Big) < \infty  
$$
 for $\, f_1, f_2, \cdots \,$ and  all  its subsequences and permutations. \qed

%%%%%%%%%%%%%%%%%%%%%%%
\subsection{Proof of Theorem \ref{thm2b} and  of Corollary \ref{cor2b}}
\label{sec4.2}
%%%%%%%%%%%%%%%%%%%%%%%

Following   the trail of \cite{BerPet}, we fix the subsequence $\, f_{k_1}, f_{k_2}, \cdots\,$ of simple functions from the construction leading   to \eqref{6.4}, and relabel it $\, f_1, f_2, \cdots\,$ for simplicity. This sequence is adapted to the   filtration $\mathbb{F}$ of \eqref{4.3a}, and a martingale difference with respect to it. We recall also  the  \textsc{Egorov}  sets $\, E_1 \subseteq  E_2 \subseteq  \cdots\,$ in $\F$ from   below \eqref{6.6};  and assume   $\,\P \big( E_k) \ge 1 - \eps_k\,,$   $\,\sup_{n \in \N} \,\E \big( f_n^2 \cdot \Ind_{\Omega \setminus E_k} \big) \le \, \eps_k\,, $   $\, \forall ~ k \in \N\,.$

In the next two subsections    we shall   construct a subsequence $\, f_{\ell_1}, f_{\ell_2}, \cdots\,$ of the relabelled $\, f_1, f_2, \cdots\,$ which is strongly $\Ll^2-$exchangeable at infinity,   with speed   $\,  \boldsymbol{\eps}=\big( \eps_k \big)_{k \in \N}\,.$

%%%%%%%%%%%%%%%%%%%%%%%
\subsubsection{The Induction Step}
\label{sec6.1.1}
%%%%%%%%%%%%%%%%%%%%%%%

We establish here Corollary \ref{cor2b}, which is a rephrasing of Theorem \ref{thm2b}. We proceed by induction on $k$.   

 Starting with $A^{(0)}_0 = \Omega$ for $k=0$, suppose that the partition $\, \big( A^{(k-1)}_j,~ j=0,1,\cdots, J_{k-1} \big)\, $ has been constructed, as well as $\, f_{\ell_1}, f_{\ell_2}, \cdots , f_{\ell_{k-1}}\,$ with the desired properties. We construct the next partition  level
  $\, \big( A^{(k)}_j,~ j=0,1,\cdots, J_{k} \big)\,$ as follows:

   Let $A^{(k)}_0:= A^{(k-1)}_0 \cap \big(\Omega \setminus E_k \big)  \,.$ Using the uniform $\, {\cal W}_2-$convergence of   $\,\big(  \boldsymbol{\mu}_{( n)}  \big)_{n \in \N}\,$ as in \eqref{6.5}-\eqref{6.6}  to the ``aggregator"   $\,\boldsymbol{\mu}: \Omega \to {\cal P}_2 (\R)\,$, whose restriction to $E_k$ has relatively compact   range, we   find an  integer $J_k > J_{k-1}\,$,   and a partition  $\, \big( A^{(k)}_j,~ j=0,1,\cdots, J_{k} \big)\,$ of $\,\Omega\,$ which \\ {\bf .} has $A^{(k)}_0$ as first element;  \\ {\bf .}   refines   the previous-level  partition $\, \big( A^{(k-1)}_j,~ j=0,1,\cdots, J_{k-1} \big)\,$ and the $\sigma-$algebra $\F_{k-1}$ of \eqref{6.7FF}; and is such that  \\ {\bf .}    {\it for every $\,j= 1,\cdots, J_{k} \,,$  the restrictions to $A^{(k)}_j$ of the measure-valued mappings $\,\big( \boldsymbol{\mu}_{( n)} \big)_{n \ge  \ell_k} \,$ and $\,\boldsymbol{\mu}\,, $ all lie in a set $\,     {\cal M}^{(k)}_j \subset {\cal P}_2 (\R)\,$ of $\,{\cal W}_2-$diameter less than $\sqrt{\eps_k\,}$.  } 
  
  \smallskip
   Furthermore,   we may (and do) assume  the functions $\, f_{\ell_1}, f_{\ell_2}, \cdots , f_{\ell_{k-1}}\,$ to be measurable with respect to the $\sigma-$algebra generated by the partition  $\, \big( A^{(k)}_j,~ j=0,1,\cdots, J_{k} \big)$.          Consequently, and in the notation of \eqref{6.2}, \eqref{6.4}, we may find an integer $\, \ell_k > \ell_{k-1}\,$ with the property 
  $$
  {\cal W}_2 \Big( \mu^{A^{(k)}_j}_{(n)} , \mu^{A^{(k)}_j} \Big) \,\le\, \sqrt{\eps_k\,}\,, \qquad \forall ~~ n \ge  \ell_k \,, ~ j=1, \cdots, J_k\,.
  $$
Results  of  \textsc{Berkes-Philipp} (\cite{BerPhi}, Theorems 1,\,2) along with  the assumed richness of the $\sigma-$algebra $\F$ and properties of the quadratic \textsc{Wasserstein} distance (the ``joining step" in   \S\,\ref{sec6.1.2a}),  provide now  a sequence  of independent functions $ \,h^{A^{(k)}_j}_{k+1},  \,h^{A^{(k)}_j}_{k+2}\,,$ with common distribution   $\, \mu^{A^{(k)}_j}$  as in   \eqref{6.4}   under  the probability measure  $\,\P^{\,A^{(k)}_j}   $ of  \eqref{6.7}, and  the property    $$ \E^{\,\P^{\,A^{(k)}_j}} \Big( h^{A^{(k)}_j}_{n}-f_{\ell_n} \Big)^2 \le \,\varepsilon_k\,,  \quad  \forall ~ ~   n \ge  \ell_k \,, ~   ~  j=1, \cdots, J_k\,. $$

%%%%%%%%%%%%%%%%%%%%%%%
\subsubsection{The Joining  Step}
\label{sec6.1.2a}
%%%%%%%%%%%%%%%%%%%%%%%

The above arguments need the following  simple property of the \textsc{Wasserstein} distance.  

\smallskip
\noindent
{\it Suppose two measures  $\mu\,, \,\nu$ in ${\cal P}_2 (\R)$ satisfy $\, {\cal W}_2^2 \big( \mu, \nu \big) <  \eps \,,$ and that a given function $f \in \Ll^2 (\Omega, {\cal F}, \P)$ has distribution $\mu$. With $\,\big(\,\overline{\Omega}, \overline{{\cal F}}, \overline{\P}\,\big) := \big( [0,1] \times \Omega, \,{\cal B} ([0,1]) \otimes{\cal F}, \,\mathrm{Leb} \otimes \P)\,,$ there exists a function $ \, g \in \Ll^2 \big(\,\overline{\Omega}, \overline{{\cal F}}, \overline{\P}\,\big)$ with distribution $\nu\,,$  and   such that $\, \overline{\E} \,\big( f-g \big)^2 < \eps\,.$  
}

 \smallskip 
The verification of this property is particularly  straightforward, when the function $f $ is simple. As we only need this case,  we take   $\, f = \sum_{j=1}^N \alpha_j \, \Ind_{A_j}\,$ for some real numbers $\, \alpha_1, \alpha_2, \cdots , \alpha_N\,$ and a finite partition $\, A_1, A_2, \cdots, A_N\,$ of $\Omega\,.$ By definition  of the  \textsc{Wasserstein} distance, there are probability measures $\, \kappa_1, \kappa_2, \cdots, \kappa_N\,$ in ${\cal P}_2 (\R)$ with 
$  
\,\nu = \sum_{j=1}^N  \, \P \big( A_j\big)\, \kappa_j\,,  ~ ~\sum_{j=1}^N  \, \P \big( A_j\big)\, {\cal W}_2^2 \big( \boldsymbol{\delta}_{\alpha_j} \,,\kappa_j\big) < \, \eps\,.
 $

Now, for every $\, j = 1, \cdots, N,$ there exists a function $\, g_j \in \Ll^2 \big( [0,1], {\cal B} ([0,1]) , \,\mathrm{Leb} \big)\,$ with distribution $\, \kappa_j\,;$    we use these   to define   
$\,
g(t, \omega) := \sum_{j=1}^N \, g_j (t) \,  \Ind_{A_j} (\omega)\,,~  (t, \omega) \in [0,1] \times \Omega\,.
$ This new function $g$ has the desired property $\, \overline{\E} \,\big( f-g \big)^2 < \eps\,.$  \qed

%%%%%%%%%%%%%%%%%%%%%%%
\section{The Proof of Theorem \ref{H-HRE}\,(i)}
\label{sec5}
%%%%%%%%%%%%%%%%%%%%%%%

We recall subsection 2.3 and the first few paragraphs of section \ref{sec6}. On their strength,    it is enough to consider  functions   $\, f_1, f_2, \cdots$  whose squares are uniformly integrable;  and   reasoning as in the preamble of section \ref{sec4} (cf.\,section \ref{sec7}), assume these are  simple, centered martingale differences with $f_\infty \equiv 0$.  We appeal now  to Theorem \ref{prop2b},     recalling its bounded-in-$\Ll^2\,$ approximating  and  omnibus sequences $f_{1}^*, f_{2}^*, \cdots $  and $h_1, h_2, \cdots$, respectively, and  to subsection \ref{rem2.4}; we  introduce also the functions
\begin{equation} 
\label{5.0}
  \xi_n \, :=\, f^*_n - h_n\,, \quad n \in \N\,.
\end{equation}
Denoting   by $\boldsymbol{\zeta}  \in \Ll^1$ the weak-$\Ll^1$ limit   of    $\, \xi_1^2, \xi_2^2, \cdots$  in the manner of \eqref{2.2}, we   observe from \eqref{6.8}  the bound $\, \E \big[ \,\xi_n^2 \, \big|\, {\cal B} \,\big] \le 1\,,$ valid for every $ n   \in \N$  and leading to     $\,\P \big( 0 \le \boldsymbol{\zeta} \le 1\big) =1\,.$      
 
 \smallskip
 {\it At  this point, Proposition \ref{H-HRE-1} takes over:} after passing once again to a (relabelled) subsequence, Proposition \ref{H-HRE-1} applies to the sequence $\xi_1, \xi_2, \cdots$  (and to all its subsequences and permutations)   and gives 
\begin{equation} 
\label{5.1}
\sum_{N \in \N} \, \P \left( \bigg| \frac{1}{N} \sum_{n=1}^N \xi_n   \bigg| > \frac{\,\eps\,}{4}   \right)< \infty\,, \quad \forall ~ \eps > 0\,.
\end{equation}
The omnibus  sequence $h_1, h_2, \cdots$ satisfies the tenets of Proposition \ref{prop3}, as do all its subsequences and permutations, thus also its conclusion  
$ \, \sum_{N \in \N} \, \P \big( \big|   \sum_{n=1}^N h_n   \big| >  \eps\,N / 4   \big)< \infty\,$,   which leads to 
$$
\sum_{N \in \N} \, \P \left( \bigg| \frac{1}{N} \sum_{n=1}^N f^*_n   \bigg| > \frac{\,\eps\,}{2}   \right)< \infty \,,\qquad \forall~  \eps >0
$$
  on account of  \eqref{5.0}--\eqref{5.1}. Therefore, and in   the context of Theorem \ref{prop2b} again,  in order to establish the HRE property $\, \sum_{N \in \N} \, \P \, \big( \big|   \sum_{n=1}^N f_{\ell_n}   \big| >  \eps\,N    \big) < \infty\,, ~ \eps >0\,$ (for the subsequence $f_{\ell_n}, f_{\ell_n}, \cdots$) it suffices to show $\, \sum_{N \in \N} \, \P \left(  \,  \sum_{n \in \N}  \big| f_{\ell_n}  - f^*_n   \big| >  \eps  N  / 2 \right) < \infty\,.$ 
  
  \smallskip
  But the elementary ``layered representation of the expectation" 
\begin{equation} 
 \label{4.31}
\sum_{N \in \N} \P (Z > N) \,\le \,\E (Z)\, = \int_0^\infty \P \big( Z > t \big)   \ud t \, \le \,   \sum_{N \in \N_0} \P (Z > N)\,,\quad \forall ~ Z \in \Ll^0_+\,,
\end{equation}  
leads to the   bounds  
$$
\frac{\eps}{\,2\,} \cdot   \sum_{N \in \N} \, \P \left(  \,  \sum_{n \in \N}  \Big| f_{\ell_n}  - f^*_n   \Big| > \frac{\,\eps\,}{2}\,N   \right)  \le \,     \E   \,\sum_{n \in \N}  \Big| f_{\ell_n}  - f^*_n   \Big| \,  \le \,    \sum_{n \in \N}  \left( \E \, \Big| f_{\ell_n}  - f^*_n   \Big|^2 \right)^{1/2}  < \infty
$$
  on account of   $\, \E \,\big( f_{\ell_n}- f_{ n}^*  \big)^2 \le 2^{-n}\,$ from Theorem \ref{prop2b}, and completes the argument for the proof of Theorem \ref{H-HRE}\,(i). \qed

 %%%%%%%%%%%%%%%%%%%%%%
\section{The Proof of Theorem \ref{H-HRE}\,(ii)}
\label{sec3d}
%%%%%%%%%%%%%%%%%%%%%%%

We   turn  now to Part (ii) of Theorem \ref{H-HRE}. Namely, we consider a sequence $\, f_1, f_2, \cdots\,$ of measurable, real-valued functions   satisfying the HRE property \eqref{2.4} for some  $\, f_\infty \in \Ll^2\,$,  i.e., 
 $$
   \sum_{N\in \N} \,\P\, \bigg( \bigg| \sum_{n=1}^N f_{k_n}  - N \, f_\infty \bigg| > \eps \, N \bigg)  \, < \, \infty \,, \quad \forall ~ \eps >0    
 $$
along some   subsequence    $\, f_{k_1} ,f_{k_2},\cdots\,$   and    all   its subsequences and permutations.  Here, we may assume    $\,f_\infty \equiv 0\,$.

 Indeed,  $\,f_{k_1} -  f_\infty, \,f_{k_2} -  f_\infty, \cdots$ satisfies the  property \eqref{2.4} with limiting function identically equal to zero. Clearly, the sequence $\,\big( f_{k_1}-  f_\infty \big)\, \Ind_{A_1} , \,\big( f_{k_2}-  f_\infty \big) \, \Ind_{A_2}, \cdots\,$ is bounded in $\Ll^2$ if, and only if, $\,f_{k_1}\, \Ind_{A_1}    , \,f_{k_2} \, \Ind_{A_2} , \cdots$ is; and likewise,   
$\,\big( f_{k_1}-  f_\infty \big) \, \Ind_{A_1^c}\,, \big( f_{k_2}-  f_\infty \big)\, \Ind_{A_2^c} , \cdots \,$ converges to zero in $\Ll^1$ if, and only if, $\,f_{k_1}\, \Ind_{A_1^c}    , \,f_{k_2} \, \Ind_{A_2^c} , \cdots \,$ does.

Because complete convergence  implies  convergence a.e., this  leads to the a.e.\,convergence of the \textsc{Ces\`aro} averages in \eqref{1.8}   to $f_\infty \equiv 0\,,$ also 
hereditarily.  Whereas, passing to $\, f_{k_1} / \eps ,f_{k_2}/ \eps, $ $\cdots\,,$ it is enough to require  these hereditary  properties  only for $\,\eps =1$, namely,   
\begin{equation} 
\label{2.9_four}
  \sum_{N\in \N} \,\P\, \bigg( \bigg| \sum_{n=1}^N f_{k_n}  \bigg| >   N \bigg)  \, < \, \infty \,.   
 \end{equation}
 We shall  show   that   there exist       sets $\, A_1 \,,\, A_2 \,, \cdots\,$ in $\,  {\cal F} \,$    with    $\, \lim_{n \to \infty} \,\P (A_n) = 1\, $        (cf.\,\eqref{Z.26} below),  as well as  a further, relabelled  subsequence $ f_{k_1}, f_{k_2}, \cdots\,, $  such that          $  f_{k_1} \mathbf{ 1}_{A_{k_1}}, f_{k_2} \mathbf{ 1}_{A_{k_2}} \cdots   $  is bounded in $\, \Ll^2  ,$  while   $  f_{ k_1} \mathbf{ 1}_{A_{k_1}^c}, f_{ k_2} \mathbf{ 1}_{A_{k_2}^c}, \cdots  $ converges to zero  in $\,\Ll^1.$

 %%%%%%%%%%%%%%%%%%%%%%%
\subsection{The Plan}
 \label{sec6.1}
%%%%%%%%%%%%%%%%%%%%%%%

This program  will be carried out in four distinct steps.

\smallskip
\noindent
{\bf Step 1:} {\it    The $  f_{k_1}, f_{k_2}, \cdots $ may be taken  bounded in $\, \Ll^0\,:$   $\, \lim_{\lambda \to \infty} \sup_{n \in \N} \P ( | f_{k_n}| > \lambda ) \to 0\,.$ }

 \smallskip
   \noindent
   {\bf Step 2:} {\it The $\, f_{k_1}, f_{k_2}, \cdots\,$ may be assumed  integrable.}

 \smallskip
\noindent
{\bf Step 3:} {\it The function $\, \boldsymbol{\eta} : \Omega \to [0, \infty]\,$ of \eqref{2.8} satisfies $\, \E (  \boldsymbol{\eta} ) < \infty\,.$} 

 \smallskip
\noindent
{\bf Step 4:} \,{\it Steps 1--3 lead   to the claims in Part (ii) of Theorem \ref{H-HRE}.}

  %%%%%%%%%%%%%%%%%%%%%%%
\subsubsection{Step 1}
 \label{sec6.1.1}
%%%%%%%%%%%%%%%%%%%%%%%

Assume the $\, f_{k_1}, f_{k_2}, \cdots\,$ were not bounded in $\, \Ll^0\,$. Then for some constant $\alpha >0$ and (relabelled) subsequence  
we would have $\,\P \big( \big| f_{k_N} \big| > \, 2\,N \big) \ge \alpha \,, ~ \forall ~ N \in \N\,;\,$ and reasoning as in \textsc{Erd\H{o}s} \cite{Erd},  also   
\begin{equation} 
\label{Z.1a}
   \Big\{ \big| f_{k_N} \big| > 2\, N \Big\} \, \subseteq \, \left\{ \bigg| \sum_{n=1}^N f_{k_n} \bigg| > N \right\} \cup \left\{ \bigg| \sum_{n=1}^{N-1} f_{k_n} \bigg| > N-1 \right\}\,, \qquad \forall ~ N \ge 2\,.
\end{equation} 
 This  then implies $\, \sum_{N \in \N } \P \big(      \big| \sum_{n=1}^N f_{k_n} \big| > N \big) = \infty \,,$ contradicting the  HRE property.

%%%%%%%%%%%%%%%%%%%%%%%%
   \subsubsection{Step 2}
    \label{sec6.1.2}
%%%%%%%%%%%%%%%%%%%%%%%

The HRE property \eqref{2.9_four} for  $\, f_{k_1}, f_{k_2}, \cdots\,,$ gives
$$
  \sum_{N\in \N} \,\P\, \bigg(   \sum_{n=1}^N f_{k_n}   >   N \bigg)  \, < \, \infty \,, \qquad  \sum_{N\in \N\,,\, N \ge 2} \,\P\, \bigg(   \sum_{n=2}^N f_{k_n}   >   N \bigg)  \, < \, \infty \,,
$$
thus also $\, \sum_{N\in \N\,,\, N \ge 2} \,\P\, \big(    f_{k_1}   >  2\, N \big) < \infty \,.$ Applying the same reasoning to $\, -f_{k_1}, -f_{k_2}, \cdots,$ and adding, we are  led  to $\, \sum_{N\in \N } \,\P\, \big(  \big|  f_{k_1}  \big| >  4\, N \big) < \infty \,,$   thus also via the "layered representation"  \eqref{4.31} to $\, \E \big(  \big|  f_{k_1}  \big|  \big) < \infty \,.$ Similar arguments lead 
to $\, \E \big(  \big|  f_{k_n}  \big|  \big) < \infty \,$ for all $\,n \in \N\,.$

 %%%%%%%%%%%%%%%%%%%%%%%
\subsubsection{Step 3}
 \label{sec6.1.4}
%%%%%%%%%%%%%%%%%%%%%%%

On the strength of Step 1, and after passing to an appropriate subsequence, we may assume that the     functions 
$f_1, f_2, \cdots$ are "determining",  in the sense of    satisfying       for $\P-$a.e.\,$\,\omega \in \Omega\,$ the stable convergence     \eqref{2.7}   for some limiting probability distribution function $  \mathbf{ H}  (\cdot\,, \omega)\,$ and corresponding random  measure $\, \boldsymbol{\mu} (\omega)\,$. As in \eqref{2.8}, we let   $\boldsymbol{\eta} (\omega)= \int_\R  \, x^2\, \ud  \mathbf{ H}  (x, \omega) \le \infty\,.$

We  consider at this point a sequence $g_1, g_2, \cdots\,$ of random variables, conditionally independent and with common distribution $\boldsymbol{\mu}\,,$ given the sigma-algebra $\, \sigma (\boldsymbol{\mu})\,$ (cf.\,\cite{A}, p.\,72 for the requisite construction). We recall also the following important complement to \eqref{2.7}, due  
to \textsc{Dacunha-Castelle} \cite{DC} (cf.\,\cite{A},      Corollary 8;\,\cite{A1},\,p.\,122): \,{\it For every $\,M \in \N\,,$ and    passing again to a suitable subsequence of $ \,f_{ 1}, f_{ 2} , \cdots \,,  $  the random vector}
\begin{equation} 
\label{Z.5}
\big( f_{k+1}, \cdots,  f_{k+M} \big) \quad \text{{\it converges in distribution as}} ~\,k \to \infty\, ~ \text{{\it to}} \quad \big( g_1, \cdots, g_M \big).
\end{equation}
 $\bullet~$ For expository reasons, we consider  a special case first. Assume the following sharpened version of the   HRE property \eqref{2.9_four}: namely,    that   there exists a constant $\, K \in (0, \infty) \,$ with  
  \begin{equation} 
 \label{Q.1}
  \sum_{N\in \N} \,\P \Big( \Big| \sum_{n=1}^N f_{k_n}  \Big| >   N \Big)   \,\le \, K  
 \end{equation} 
valid  for every subsequence $\, f_{k_1}, f_{k_2}, \cdots\,  $ of $ \,f_{ 1}, f_{ 2} , \cdots \,.  $  With  $\,R \in \N_0\, $ fixed, we obtain then 
$$
\sum_{N\in \N} \,\P\, \bigg( \bigg| \sum_{n=1}^N f_{ R+n}  \bigg| >   N \bigg)   \le  K    \,; \qquad \text{thus also} \qquad \sum_{N=1}^M \,\P\, \bigg( \bigg| \sum_{n=1}^N f_{ R+n}  \bigg| >   N \bigg)   \le  K   \,,
$$
for every given $M \in \N\,.$   Letting $\, R \to \infty\,$ and recalling \eqref{Z.5}, we obtain 
$$
\sum_{N=1}^M \,\P\, \bigg( \bigg| \sum_{n=1}^N g_{n }  \bigg| >   N \bigg)   \le  K  \,, ~~~ \forall ~ M \in \N\,, \qquad \text{thus} \qquad \sum_{N \in \N}  \,\P\, \bigg( \bigg| \sum_{n=1}^N g_{n }  \bigg| >   N \bigg)  \le K \,.
$$
But  when conditioned on the sigma-algebra $  \sigma (\boldsymbol{\mu}) ,$ the functions $\,g_1, \,g_2, \cdots\,$  are independent     with common random distribution function $ \mathbf{ H}  (\cdot\,, \omega)\,.$ For this sequence, applying     \eqref{8.000} of Proposition  \ref{prop8.0} conditionally on $\, \sigma (\boldsymbol{\mu})\,$ and taking expectations, leads to 
 $\,K \ge \sum_{N \in \N}  \,\P\, \big( \big| \sum_{n=1}^N g_{n }  \big| >   N \big) \ge   c  \cdot \E \big( g_1^2 \big)
 $   with   a  suitable universal     constant  $\,c >0\,.$   We have thus shown that, under \eqref{Q.1}, the ``randomized second moment"   $\,\omega \mapsto \boldsymbol{\eta} (\omega)= \int_\R   x^2\, \ud  \mathbf{ H}  (x, \omega) $ $\in [0, \infty]\,$ in \eqref{2.8} is   integrable, namely, satisfies    $\,  \E (\boldsymbol{\eta}) = \E(g_1^2) < \infty\,. $  
 
A similar argument  establishes  Proposition \ref{H-HRE-Three} as well.

\smallskip
\noindent
$\bullet~$ We   drop  now the condition  \eqref{Q.1}; assume only  
\eqref{2.9_four}, i.e., $ \sum_{N\in \N}  \P  \big( \big| \sum_{n=1}^N f_{k_n}  \big| >   N \big)    <   \infty \,,$ for every subsequence    $\, f_{k_1} ,f_{k_2},\cdots$ of      $\, f_{ 1} ,f_{ 2},\cdots\,;$ and  show that this leads   to $\,\E (\boldsymbol{\eta}) \, <  \, \infty\,$.  

We    argue this  by   contradiction: namely, assume   $\,\E (\boldsymbol{\eta}) = \infty\,,$         and work toward finding a subsequence    $\, f_{\ell_1} ,f_{\ell_2},\cdots$ of      $\, f_{ 1} ,f_{ 2},\cdots\,$ and an increasing  sequence of integers $1 < M_1 < M_2 < \cdots\,$   with 
\begin{equation} 
\label{Z.12}
 \sum_{N=1}^{M_j} \,\P\, \bigg( \bigg| \sum_{n=1}^N f_{\ell_n}  \bigg| >   N \bigg)  \, \ge \, j-1 \,, \qquad \forall ~ j  \in \N\,; 
\end{equation}
this will then contradict  the assumed   HRE property \eqref{2.9_four} for the sequence $\, f_{\ell_1} ,f_{\ell_2},\cdots$.

As a first reduction step, we may {\it assume that each $f_n$ is bounded,} i.e., $f_n \in \Ll^\infty.$ Indeed, as each $f_n$ can  be assumed integrable on the strength of Step 2, we may find $f_n^* \in \Ll^\infty $ with $\, \big\| f_n - f^*_n \big\|_1 < 2^{-n}.$ Now   $\, f_{ 1}^* ,f_{ 2}^*,\cdots\,$ inherits the HRE property from $\, f_{ 1} ,f_{ 2},\cdots\,$,  as  $\big\{f_n - f_{ n }^* \big\}_{n \in \N}\,$   satisfies this property on the strength of Proposition 7.2\,(i); and $\, f_{ 1}^* ,f_{ 2}^*,\cdots\,$ is still determining, with the same exchangeable sequence  $\, g_{ 1} ,g_{ 2},\cdots\,$ as the original  sequence $\, f_{ 1} ,f_{ 2},\cdots\,$. In conclusion, we shall assume that $\, f_{ 1} ,f_{ 2},\cdots\,$ are all in $  \Ll^\infty$.

 \smallskip
 Next, we return   to the task of finding a subsequence    $\, f_{\ell_1} ,f_{\ell_2},\cdots$ of      $\, f_{ 1} ,f_{ 2},\cdots\,$ with the  properties spelled out in the previous two paragraphs; {\it and proceed to prove \eqref{Z.12} via  an induction} (on $j \in \N$). 
 
 The first step $j=1$ of this induction is clear. We then assume   that \eqref{Z.12} holds for some $j \in \N\,$ and natural numbers $\,M_1 <M_2   < \cdots < M_j\,,$   $\,\ell_{1}    <   \ell_{2} < \cdots < \ell_{M_j}\, .$    For the function    $\, F_j := \sum_{n=1}^{M_j} f_{\ell_n} \in \Ll^\infty\,,$ we can find   a real constant $C_j > 0$ so large,   that $\, \big| F_j \big| \le C_j \,$ holds a.e. Summoning   the sequence $\, g_{ 1} ,g_{ 2},\cdots\,$ of random variables which are conditionally independent and with common distribution $\boldsymbol{\mu}\,$ given the sigma-algebra $\, \sigma (\boldsymbol{\mu})\,,$ applying Lemma \ref{prop8.1} conditionally on $\, \sigma (\boldsymbol{\mu}),$ taking expectations, and using monotone convergence as well as the standing assumption $\,\E (\boldsymbol{\eta}) = \infty\,,$ we obtain 
$$
   \sum_{ N \in \N \atop N \ge  M_j +1     } \,\P \bigg(\, \bigg| \,  \sum_{n=M_j+1}^N g_n \,\bigg| >    N + C_j  \bigg)  \,\ge\, c \cdot    \E (\boldsymbol{\eta}) - 1 = \infty\,.
   $$
We  extend now the string $\,\big( \ell_{1}   ,  \cdots ,\ell_{M_j} \big)\, $ to the increasing string  $\,\big( \ell_{1}   ,  \cdots , \ell_{M_j},   \ell^{R}_{M_j +1}, \cdots, $ $\ell^{R}_{M_j +K} \big) :=\big( \ell_{1}   ,  \cdots ,\ell_{M_j},$ $   R+1 , \cdots, R+K  \big) \,, $ with $\, \ell_{M_j} \le R\,$ and  $\,K \in \N\,$  fixed, so that   the     random vector $ \big(    f_{\ell^{R}_{M_j +1}}, \cdots, f_{\ell^{R}_{M_j +K}} \big)$ converges in distribution,  as   $    R \to \infty\,, $ to the      random vector $ \big(    g_{M_j +1}, \cdots, g_{M_j +K} \big)$. This  yields   
  $$
 \lim_{K \to \infty}   \lim_{R \to \infty} \sum_{N=M_j+1}^{M_j+K} \,\P\, \bigg( \bigg| \sum_{n=1}^{M_j} f_{\ell_n} + \sum_{n=M_j+1}^{M_j+K}  f_{\ell_n^{R}}  \bigg| >   N \bigg)~~~~~~~~~~~~~~~~~~~~~~~~~~~~~~~~~~~~~~~~~~~
 $$
 $$
~ \ge   \lim_{K \to \infty}   \lim_{R \to \infty}   \,\sum_{N=M_j+1}^{M_j+K} \,\P\, \bigg( \, \bigg|  
 \sum_{n=M_j+1}^{M_j+K}  f_{\ell_n^{R}}  \bigg| >   N + \,C_j  \bigg)  
 $$
 $$
  =   \lim_{K \to \infty}   \, \sum_{N=M_j+1}^{M_j+K} \,\P\, \bigg( \bigg| \sum_{n=M_j+1}^{M_j+K} g_n  \bigg| >   N+ C_j     \bigg)    = \infty\,.
 $$
Taking $R$ and $M_{j+1}$ large enough, we define 
$ \big(    \ell_{M_j +1}, \cdots, \ell_{M_{j +1}} \big) :=   \big(    \ell^{R}_{M_j +1}, \cdots, \ell^{R}_{M_{j +1}} \big)$    to obtain 
 $\,
  \sum_{N = M_j +1}^{ M_{j+1}} \,\P\, \big( \big| \sum_{n=1}^N f_{\ell_n}  \big| >   N \big)  \, \ge \, 1\,.
 $
 
 \smallskip
This is    the inductive step $\, j \mapsto j+1\,$ needed for the claim \eqref{Z.12}, and   completes    Step 3.

  %%%%%%%%%%%%%%%%%%%%%%%
\subsubsection{Step 4}
 \label{sec6.1.5}
%%%%%%%%%%%%%%%%%%%%%%%

From  Step  1, we may assume   the  sequence $f_1, f_2,   \cdots$ to be determining; in this manner,     
$\,
\lim_{n \to \infty} \P \big( f_n \le x \big) = H(x) := \int_\Omega  \mathbf{ H}  (x, \omega) \, \P (\mathrm{d} \omega) \,, ~ \forall ~ \,x \in  \mathbf{ D }
\,$ holds from \eqref{2.7} with $B = \R\,$. Whereas, Step 3   gives
 \begin{equation} 
\label{Z.53}
\lim_{n \to \infty} \E \Big( f_n^2 \cdot \Ind_{ \{ |f_n |  \le K \} } \Big) \, =  \, \int_{[-K,K]} x^2 \, \mathrm{d}  H(x) 
  \end{equation}  
  $$
 ~~~~~~~~~~~ ~~~~~~~~~~~~~
  \,=\, \E \int_{[-K,K]} x^2 \, \mathrm{d}   \mathbf{ H}  (x, \omega) \, \le \, \E \big(\boldsymbol{\eta} \big)\,=: \,M < \infty
$$
provided $\, \pm \, K   \in  \mathbf{D} \,.$   We select now a sequence $\, 0 < K_1 < K_2 < \cdots\,$ with $\, \pm \, K_n \in  \mathbf{D}\,$ and  
 \begin{equation} 
\label{Z.33}
 \P \big( |f_n| > K_n \big) \le 2^{-n}\,, \qquad \E \big( \, \big| f_n \big| \cdot \Ind_{ \{ |f_n| > K_n \} } \, \big) \le 2^{-n}\,  , \qquad \forall ~~ n \in \N
  \end{equation}  
(the latter because the $f_1, f_2, \cdots\,$ can be assumed integrable, on the strength of Step 2). Thus,     after passing again to a subsequence,    \eqref{Z.53} gives   
\begin{equation} 
\label{Z.6}
 \E \big[ \, f_n^2 \cdot \Ind_{ \{ |f_n| \le K_n \} } \, \big] \,< \,2\,M\, \quad \text{ for all $n \in \N\,.$}
 \end{equation} 
{\it We introduce at this point the sought-after sets}
  \begin{equation} 
\label{Z.26}
     A_n := \big\{  |f_n| \le K_n \big\}\, \quad \text{for} \quad     n \in \N\,,\quad \text{which satisfy} ~~ \P (A_n) > 1-2^{-n}
    \end{equation} 
    (passing to   $  \widetilde{A}_n = \bigcap_{\,m \ge n} A_m $ if necessary, we may   assume  that these     $ \, A_1, A_2, \cdots\, $     are  increasing).  From  \eqref{Z.6},   the sequence   $ f_1 \, \Ind_{A_1}, f_2 \, \Ind_{A_2}, \cdots$  is bounded in $\Ll^2,$  as posited in Theorem    \ref{H-HRE}\,(ii); whereas    $ f_1 \, \Ind_{A_1^c}, f_2 \, \Ind_{A_2^c}, \cdots\,$  converges to zero  in $\Ll^1,$  from the second display in \eqref{Z.33}.  

\smallskip
 Step 4 is thus established, and  the Proof of Theorem \ref{H-HRE}\,(ii) is   complete.  \qed

%%%%%%%%%%%%%%%%%%%%%%%
\section{Appendix: Functions  Supported on Disjoint  Sets}
\label{sec3}
%%%%%%%%%%%%%%%%%%%%%%%

We recall     from \cite{CM}, Lemma 2.1.3 (cf.\,\cite{KP};\,\cite{KK}, Lemma A.44;\,\cite{Ro}) the important  "KPR Lemma".

\begin{lemma} 
[\textsc{Kade\'c-Pe\l zy\'nski-Rosenthal}]
\label{KPR}
Every bounded--in--$\Ll^1$ collection of real-valued functions    contains  a   sequence $f_1, f_2,   \cdots$ of the form   $f_n = h_n + g_n\,,~ n \in \N\,, $ where the $h_1, h_2,   \cdots$ are  integrable functions    supported on disjoint sets $B_1, B_2, \cdots\,,$  and  the   $g_1, g_2,   \cdots$ are uniformly integrable functions.
\end{lemma}

 The following   result  shows that fast convergence to zero in $\Ll^1$, is sufficient for the HRE property; as well as necessary, after passage to an appropriate  subsequence,   for functions with disjoint or "essentially disjoint" (in a manner reminiscent of 
 Lemma \ref{KPR}) supports.  

 \begin{proposition} 
 \label{lem3.7} {\rm Conditions Sufficient,  and   Necessary, for the HRE Property.} \\
Consider a  sequence of  real-valued, measurable    functions  $\, h_1, h_2, \cdots\,$.  

\noindent
{\bf (i)} The condition  of fast convergence to zero in $\Ll^1$, namely, 
\begin{equation} 
\label{3.12}
\sum_{n \in \N} \E \big( \big| h_n \big|   \big) < \infty
\end{equation} 
implies, for $\, h_1, h_2, \cdots\,$   and   all its subsequences and permutations, the  complete convergence  
 \begin{equation} 
\label{3.5.a}
\sum_{N \in \N} \, \P \Big( \Big|   \sum_{n=1}^N  h_n  \Big| > \eps N   \Big)< \infty\,, \qquad \forall ~ \eps > 0\,.
\end{equation}
    {\bf (ii)}  Conversely, if  the  $\,h_1, h_2, \cdots$ are     supported on   disjoint sets $\,B_1, B_2, \cdots\,$ in ${\cal F}$, then, after passing to an appropriate subsequence,    the   validity of  
   \eqref{3.5.a} implies 
\begin{equation} 
\label{3.13}
\underline{\lim}_{\,n \to \infty} \,\E \big( \big| h_n \big|   \big)=0\,,
\end{equation} 
 thus also \eqref{3.12} along an appropriate subsequence and   all its  subsequences and permutations.

\noindent
{\bf (iii)}  More generally, with       $\,h_1, h_2, \cdots$ as in   (ii), consider    a sequence  of measurable functions $\, g_1, g_2, \cdots\,$  with each $\,g_n$ supported on  the set $\,\Gamma_n := \bigcup_{m > n}B_m \in {\cal F}\,$. 

Then, after passing to an appropriate subsequence,  a necessary condition  for  the sequence  $\, h_1 + g_1, \,h_2 + g_2\,, \cdots\,$ to contain a subsequence   converging completely in \textsc{Ces\`aro} mean to zero,  is again \eqref{3.13}$;$ which  leads to \eqref{3.12} along a (relabelled) subsequence and all its subsequences  and permutations.
\end{proposition}

 %%%%%%%%%%%%%%%%%%%%%%%
\subsection{The Proof of Proposition \ref{lem3.7}}
\label{sec3.2}
%%%%%%%%%%%%%%%%%%%%%%%

{\it Proof of Part (i)}\,:  This follows   from  \eqref{4.31}, which gives $\,\sum_{N \in \N} \, \P \big( \big|   \sum_{n=1}^N  h_n  \big| >   N   \big) \le \sum_{N \in \N} \, \P \big(     \sum_{n \in \N} \big|  h_n  \big| >   N   \big) \le \sum_{n \in \N } \E \big( \big| h_n \big|   \big) < \infty\,.  $ This applies also to every subsequence.

 \smallskip
 \noindent
 {\it Proof of Part (ii)}\,:  We   argue by contradiction.   Suppose that \eqref{3.13} fails; to wit, that  after passing to a subsequence,  we have $\, \E\big( \big| h_n \big| \big)> 2\,\beta\,, ~\forall~   n \in \N\, $ for some $\, \beta >0\,.$    Then,   passing   to a subsequence once again  and recalling that the $\,h_1, h_2, \cdots$ are     supported on   disjoint sets, 
 \begin{equation} 
\label{3.14} 
\E \Big(  \big|  h_n \big| \cdot \Ind_{ \{   |  h_n  | >n \} } \Big) \,>\, \beta\,, \  \qquad \forall ~~n \in \N  
\end{equation}
  may also be assumed. Now, for fixed $m \in \N$, we note
\begin{equation} 
\label{3.15} 
\sum_{N \ge m } \P \big( \big| h_m \big| > N \big) \, \ge \, \E \Big(  \big|  h_m \big| \cdot \Ind_{ \{   |  h_m  | >m \} } \Big) \ge  \sum_{N > m } \P \big( \big| h_m \big| > N \big) ~~~~~~~~~~~~~  
\end{equation}
$$
~~~~~~~~~~~~~~~~~~~~~~~~   ~
= \sum_{N \ge m } \P \big( \big| h_m \big| > N \big) - \P \big( \big| h_m \big| > m \big) \ge \sum_{N \ge m } \P \big( \big| h_m \big| > N \big) - \P \big( B_m \big)
$$
(for the first two inequalities, recall \eqref{4.31}),  as well as
\begin{equation} 
\label{3.16} 
\sum_{N \in \N } \P \Big( \Big| \sum_{n=1}^N h_n \Big| \cdot \Ind_{B_m} > N \Big) = \sum_{N \ge m } \P \big( \big| h_m \big| > N \big)  \ge \ \E \Big(  \big|  h_m \big| \cdot \Ind_{ \{   |  h_m  | >m \} } \Big) > \beta
\end{equation}
from  \eqref{3.14}, \eqref{3.15}. We sum now in \eqref{3.16} over the disjoint sets $B_m\,,$ $ m \in \N\,,$ and deduce
 \begin{equation} 
\label{3.17}
 \,\sum_{N \in \N} \, \P \bigg( \Big|   \sum_{n=1}^N  h_n  \Big| >   N   \bigg)= \infty\, ;
 \end{equation}
i.e., that \eqref{3.5.a} fails for $\eps =1\,,$  thus   for all $\eps>0\,$ by scaling. This yields the   contradiction.

 \smallskip

  \noindent
{\it Proof of Part (iii)}\,: We   argue again by contradiction:   assume $\,
\underline{\lim}_{\,n \to \infty} \,\E \big( \big| h_n \big|   \big)\,>\,0\,,$
 and show that the following analogue  of 
 \eqref{3.17}  has then to hold along some subsequence:
 \begin{equation} 
\label{3.17a}
 \sum_{N \in \N} \, \P \bigg( \Big|   \sum_{\ell =1}^N \Big(  h_{k_\ell}  + g_{k_\ell} \Big)    \bigg| >   N   \bigg) = \infty\,.
 \end{equation}
    More precisely, we  suppose again that $\, \E \big( \big| h_n \big| \big) > 2\, \beta\,, ~ \forall \, n \in \N\,$ holds for some $\,\beta >0\,;$    also, without loss of generality, that   $\, \P (B_n) < n^{-3}\,$ holds for all $n \in \N\,$.  
    
    We construct now inductively a sequence of integers $1=k_1 < k_2 < \cdots\,$ with each $\,k_n\,,\, n \ge 2\,$ equal to either $2n-1$ or $2n$,  as follows. 
 For $n=1$, we take $k_1=1$ as already mentioned, so that, with the help of \eqref{4.31},     
 $$
  \,\sum_{N \in \N} \, \P \Big( \Big|  \sum_{\ell =1}^N \Big(  h_{k_\ell}  + g_{k_\ell} \Big) \Ind_{B_1}  \Big| >   N   \Big)= \sum_{N \in \N} \, \P \big( \big| h_1 \big| > N \big) \ge \E \big( \big| h_1 \big|   \big) - \P (B_1) > 2 \beta - \P (B_1) 
 $$
holds;  and note that this   is valid also along any subsequence $\,k_1, k_2, \cdots$  with  $k_1 = 1$. 
 
  Proceeding with the induction, suppose   that $1=k_1 < k_2 < \cdots < k_n\,$ have been chosen  so that, for each $\, j=1, \cdots, n\,,$  the expression 
 $$
 \sum_{N \in \N} \, \P \bigg( \Big|   \sum_{\ell=1}^N \Big(  h_{k_\ell}  + g_{k_\ell} \Big) \Ind_{B_{2j-1}}  \Big| >   N   \bigg)= \sum_{N \in \N} \, \P \bigg( \Big|   \sum_{\ell=1}^{N \wedge n} \Big(  h_{k_\ell}  + g_{k_\ell} \Big) \Ind_{B_{2j-1}}  \Big| >   N   \bigg)  
  $$
 dominates the quantity $  \,   2 \beta - \big( 2j-1 \big) \cdot \P \big( B_{2j-1} \big). $   We  pair   the function $\,h_{2n+1}\,$ with the sum $\,\sum_{\ell =1}^n g_{k_\ell} \cdot \Ind_{B_{2n+1}}\,,$ and distinguish two possibilities: 
  
   \smallskip
  \noindent
  $\bullet~$ If 
  \begin{equation} 
\label{3.18}
\E \,\bigg( \bigg| \, h_{2n+1} + \sum_{\ell =1}^n g_{k_\ell} \cdot \Ind_{B_{2n+1}} \bigg| \bigg)> \beta 
 \end{equation}
holds,  we set $\, k_{n+1} := 2n+1\,$ and proceed to estimate
 $$
 \sum_{N \in \N} \, \P \left( \Big|   \sum_{\ell =1}^N \Big(  h_{k_\ell}  + g_{k_\ell} \Big) \Ind_{B_{2n+1}}  \bigg| >   N   \right)= \sum_{N \in \N} \, \P \left( \,\left|   \sum_{\ell=1}^{N \wedge (n+1)} \Big(  h_{k_\ell}  + g_{k_\ell} \Big) \Ind_{B_{2n+1}}  \right| >   N   \right)
 $$
 $$
 ~~~~\ge  \sum_{N \in \N} \, \P \left( \bigg| \Big( h_{2n+1} +  \sum_{\ell=1}^n    g_{k_\ell} \Big) \Ind_{B_{2n+1}}  \bigg| >   N   \right)- n \cdot \P \big( B_{2n+1} \big) \ge \beta - n \cdot \P \big( B_{2n+1} \big)\,.
 $$
  $\bullet~$ If, on the other hand, \eqref{3.18} fails,  we  have $\, \E \, \big|   \sum_{\ell=1}^n    g_{k_\ell} \cdot \Ind_{B_{2n+1}} \big| > \beta\,$ on account of the assumption $\, \E \big( \big| h_{2n+1} \big| \big) > 2\, \beta\,;$     take   $\, k_{n+1} := 2n+2\,;$ and obtain  
  $$
 \sum_{N \in \N} \, \P \left( \Big|   \sum_{\ell =1}^N \Big(  h_{k_\ell}  + g_{k_\ell} \Big) \Ind_{B_{2n+1}}  \bigg| >   N   \right)= \sum_{N \in \N} \, \P \left( \,\left|   \sum_{\ell=1}^{N \wedge (n+1)} \Big(  h_{k_\ell}  + g_{k_\ell} \Big) \Ind_{B_{2n+1}}  \right| >   N   \right)
 $$ 
 $$
 ~~~~~~~~~~~~~~~~~~~~~~~~~~~~~\ge  \sum_{N \in \N} \, \P \left( \bigg|    \sum_{\ell=1}^n    g_{k_\ell} \Ind_{B_{2n+1}}  \bigg| >   N   \right)- n \cdot \P \big( B_{2n+1} \big) \ge \beta - n \cdot \P \big( B_{2n+1} \big)\,.
 $$
 Since $\, \sum_{n \in \N} \,n \cdot \P   \big( B_{2n+1} \big) < \infty\,$ holds by assumption,  we conclude by   summing up in the above display over the sets $\, B_m \,,\, m \in \N\,,$ and arrive at the desired    \eqref{3.17a}.     \qed

%%%%%%%%%%%%%%%%%%%%%%%
\section{Appendix: Approximation by Simple Martingale Differences}
\label{sec7}
%%%%%%%%%%%%%%%%%%%%%%%

We illustrate here how to approximate   bounded-in-$\Ll^2$ sequences  of functions by simple, square-integrable   martingale differences, in the manner of  \textsc{Koml\'os} \cite{K} and   \textsc{Chatterji} (\cite{Ch}, \cite{Ch2}).

Let us consider  then  a    sequence $f_1, f_2, \cdots\,$ satisfying \eqref{2.1};   this    contains   a (relabelled) subsequence converging weakly in $\Ll^2$  to some $f_\infty \in \Ll^2\,$, as in \eqref{4.1}. We take $f_\infty=0\,$ for concreteness; and approximate each $f_n$ by a {\it simple}  function  $\,h_n \in \Ll^2\,$ with 
\begin{equation} 
\label{7.1} 
\E \big| f_n - h_n \big|^2 \, \le \, 4^{-n}, ~\, \forall ~n \in \N \,; 
\end{equation}
\begin{equation} 
\label{7.1z}
  \E \bigg( \, \sum_{n \in \N} \big| f_n - h_n \big| \bigg)\, \le \,\sum_{n \in \N} \sqrt{\,  \E \big| f_n - h_n \big|^2\,}\,\le \,\sum_{n \in \N} 2^{-n} \,= \,1\,;
\end{equation}
thus also $\,  \sum_{n \in \N} \big| f_n - h_n \big| < \infty\,,$   $\P-$a.e. 
  In addition, for every test function $\xi \in \Ll^2$  we have 
$\,
\E \big( h_n \cdot \xi \big) = \E \big( f_n \cdot \xi \big) - \E \big( \big( f_n- h_n \big)  \cdot \xi \big) \, \longrightarrow \,0\,, ~ \text{as} ~~ n \to \infty
\,$
from  \eqref{4.1},   \eqref{7.1}. 
It develops that the sequence of simple functions $h_1, h_2, \cdots$   is bounded in $\Ll^2$ (since the inequality $\| h_m \|_2 \le \| f_m \|_2 + \| f_m -h_m \|_2\le \sup_{n \in \N} \| f_n \|_2+ (1/2)< \infty$ holds for all $m \in \N$), and as such  contains a (relabelled) subsequence  converging weakly in  $\Ll^2$ to $f_\infty \equiv 0$.

 \smallskip
{\it We construct now   inductively a sequence $1=k_1 < k_2 < \cdots\,$ of integers $\,k_n \ge n\,$ such that  
 \begin{equation} 
\label{7.2}
 \boldsymbol{\vartheta}_n \,:=\, \E \big( h_{k_n} \big| \,{\cal H}_{n-1} \big)\, \qquad \text{with} \qquad  {\cal H}_{n-1} := \sigma \big( h_{k_1}, \cdots h_{k_{n-1}} \big)\,, \quad n=2,3, \cdots\,,
\end{equation}
which are    simple functions, satisfy the bound $\, \big| \, \boldsymbol{\vartheta}_n \,\big| \,\le\, 2^{-n}\,, ~ \P-\text{a.e.}$}   

\smallskip
This is done as follows: the  function $\,h_{k_1} \equiv h_1\,$ is simple, so the conditional expectation of a generic $h_n$, given $h_1$, is also simple: $\, \E \big( h_{ n}\, \big| \,h_{ 1} \big)\,=\, \sum_{j=1}^J \gamma_{j}^{(n)} \cdot \Ind_{A_j}\,. $  Here, the disjoint sets $A_1, \cdots, A_J$ form a partition of $\Omega\,$; each of them has positive measure; and for every $\,j=1, \cdots , J\,,$ we have  $\,
 \gamma_{j}^{(n)} \,:=\, \big(\P (A_j)\big)^{-1} \cdot \E \big( h_n \cdot   \Ind_{A_j} \big)  \, \longrightarrow \,0\,, ~~ \text{as}  ~ n \to \infty\,. $  
 
 We can select, therefore,  $k_2>1=k_1$ so that $ \, \big|  \gamma_{j}^{(k_2)} \big| \le 2^{-2}\,, ~~ j=1, \cdots, J\,$ holds;    thus  also  $ \big|  \boldsymbol{\vartheta}_2 \big|   =  \big| \,\E \big( h_{k_2} \big| \,h_{k_1} \big)\, \big| \, =\,\sum_{j=1}^J \, \big| \gamma_{j}^{(k_2)}\big| \cdot \Ind_{A_j}\,\le\, 2^{-2}\,, ~   \P-\text{a.e.}  $  \,  Now, we    keep repeating this procedure; at each of its stages, the   vector of simple functions $\, \big( h_{k_1}, \cdots, h_{k_{n-1}} \big)$ generates a finite partition of $\Omega\,,$ and we arrive     inductively at  the claim  $\,\big| \,\boldsymbol{\vartheta}_n \,\big| \,\le\, 2^{-n}\,, ~ \P-$a.e.  \qed

 \medskip
We have now the following result.

 \begin{proposition} 
\label{prop8.2}
 For  the integers $1=k_1 < k_2 < \cdots\,$ with $\,k_n \ge n\,$ as    above, the   subsequence 
$
 \,  f_{k_1}, \,f_{k_2}, \cdots\,
$ 
of the original bounded\,--\,in\,--$\,\Ll^2$ sequence $f_1, f_2, \cdots$ has the HRE property \eqref{2.4} with $f_\infty \equiv 0$ if, and only if, this is true for the simple  martingale-differences  
\begin{equation} 
\label{7.6}
\boldsymbol{\beta}_n\,= \, h_{k_n}- \boldsymbol{\vartheta}_n\,=\, h_{k_n}- \E \big( h_{k_n}\, \big| \,h_{k_1}, \cdots h_{k_{n-1}} \big)\,,\quad n \in \N\,,
\end{equation}
which are  generated by the simple   functions $
 \,  h_{k_1}, \,h_{k_2}, \cdots\,
$
  in  \eqref{7.1}--\eqref{7.2} and    satisfy
 \begin{equation} 
\label{7.7}
\big\| \,f_{k_n} - \boldsymbol{\beta}_n \big\|_2 \, \le \, 2^{\,1-n}\,, \quad \forall ~~n \in \N\,. 
\end{equation}
 \end{proposition}

\noindent
{\it Proof:}   For the "if\," part, we note  the inequality
$$
\sum_{N \in \N}   \P \left( \bigg|  \sum_{n=1}^N  f_{k_n}  \bigg| >   N   \right)  \le  \sum_{N \in \N}   \P \left( \bigg|  \sum_{n=1}^N \Big( h_{k_n}- \boldsymbol{\vartheta}_n \Big) \bigg| > \frac{N}{\,3\,}    \right)  + \sum_{N \in \N}  \Big( \P \big( \Phi_N \big) +     \P \big( \Theta_N \big)\Big)
$$
with 
$$
\Phi_N \,:=\, \left\{ \,\bigg| \, \sum_{n=1}^N \Big( f_{k_n} - h_{k_n} \Big) \, \bigg| > \frac{N}{\,3\,}\,   \right\}\,, \qquad \Theta_N \,:=\, \left\{  \,\bigg| \, \sum_{n=1}^N \boldsymbol{\vartheta}_n\, \bigg| > \frac{N}{\,3\,}  \,  \right\}\,. 
$$
\noindent
We recall  also the   bound \eqref{4.31}, valid  for  measurable $Z \ge 0\,$; in conjunction with  $\,\big| \,\boldsymbol{\vartheta}_n \,\big| \,\le\, 2^{-n}\,$  and the resulting $ \big| \, \sum_{n \in \N} \boldsymbol{\vartheta}_n\, \big| \le 1 ,$ as well as   \eqref{7.1}, this leads to    $\, \sum_{N \in \N} \P (\Theta_N) \le 4\,$ and  $\, \sum_{N \in \N} \P (\Phi_N) \le 4\,.$  On account of the above inequality,    the "if\,"  part of the claim is thus established;    whereas,   repetition of the same argument establishes the "only if\," part, via   
$$
 \sum_{N \in \N} \P \left( \bigg|  \sum_{n=1}^N 
\Big(   h_{k_n} - \boldsymbol{\vartheta}_n \Big)   
\bigg| >   N   \right)  \le  \sum_{N \in \N}   \P \left( \bigg|  \sum_{n=1}^N   f_{k_n}  \bigg| > \frac{N}{\,3\,}    \right)  + \sum_{N \in \N}    \Big( \P \big( \Phi_N \big) +     \P \big( \Theta_N \big)\Big).  
$$
 As for  \eqref{7.7}, this follows from \eqref{7.1}, the bound $\, \big| \,\boldsymbol{\vartheta}_n \,\big| \,\le\, 2^{-n}\,,$  and the triangle inequality.  \qed

    %%%%%%%%%%%%%%%%%%%%%%%
\section{Appendix: Quantitative Results of \textsc{Hsu-Robbins-Erd\H{o}s}\,Type}
\label{sec8}
%%%%%%%%%%%%%%%%%%%%%%%

The proof of Theorem \ref{H-HRE} uses      the following uniform sharpening of the  \textsc{Hsu-Robbins-Erd\H{o}s} theorem. This   extends   a  result of \textsc{Heyde} \cite{H}, on which the characterization \eqref{1.6} is based. 

\begin{proposition} 
\label{prop8.0}
Let $\ f_1, f_2, \cdots$ be a sequence of independent copies of a random variable $f \in \Ll^1$ with $\,\E (f) =0\,,   \,\sigma^2 := \E (f^2) \le \infty\, $. For some     constants $\,0  < c < C_1< \infty\,,$ $\,0    < C_2< \infty\,$ which are {\rm universal} $($i.e.,  do not depend on the distribution of $f)$, we have then 
\begin{equation} 
\label{8.00}
c \cdot \sigma^2 \,\le  \, \sum_{N \in \N } \,\P \Big(\, \Big| \sum_{n=1}^N f_n \Big| >   N \Big) + 1 \,\le\, C_1 \cdot   \sigma^2 + C_2\,,
\end{equation} 
as well as the scaled version 
\begin{equation} 
\label{8.000}
c \cdot \sigma^2 \,\le 
 \eps^2 \bigg( \sum_{N \in \N} \,\P \bigg(\, \bigg| \sum_{n=1}^N f_n \bigg| >  \eps \, N \bigg) +1 \bigg)\,\le\, C_1 \cdot   \sigma^2 + C_2 \cdot \eps^2\,, ~~~~~~ \forall ~ \eps \in (0,1]\,.
\end{equation}
\end{proposition} 
 
\noindent
{\it Proof:} The proof of the upper bound in \eqref{8.00} follows from the inequality (47) in \textsc{Fuk-Nagaev} \cite{FN} (see also \cite{H}, p.\,175), which  gives
\begin{equation} 
\label{FN}
 \P \bigg(\, \bigg| \sum_{n=1}^N f_n \bigg| >   N \bigg)\,\le\, N \, \, \P \Big( \big|f \big| > \frac{N}{\,4\,}   \Big)+ \frac{\, 128 \, \big( 1 + 2 e^4\big)\,}{N^2  }\,, \qquad \forall   ~~ N \in \N\,;
\end{equation} 
adding over $N \in \N\,,$ we obtain the upper bound in \eqref{8.00} for some universal constants $\,C_1>0\,$, $C_2>0\,$. The lower bound     follows from the next result;    the scaled version \eqref{8.000} is obvious.   \qed

\begin{lemma} 
\label{prop8.1}
There is a universal constant $\, c \in (0, \infty)\,$ such that, for 
   independent copies $\, f_1, f_2, \cdots$  of a random variable $f  $ with arbitrary  distribution which is symmetric $($i.e.,  $-f$ has the same distribution as $f)$, and with $\, \sigma^2 := \E (f^2) \le \infty\, ,$ we have 
\begin{equation} 
\label{8.1}
   \sum_{N \in \N    } \,\P \Big(\,  \sum_{n=1}^N f_n   >2\, N   \Big) + 1\,\ge\, c \cdot   \sigma^2  \,.
\end{equation} 
Whereas, if we drop the symmetry assumption on $f$, we still have for a (possibly different)   universal constant $\, c \in (0, \infty), $ the bound
\begin{equation} 
\label{8.2}
   \sum_{N \in \N  } \,\P \bigg(\, \bigg| \sum_{n=1}^N f_n \bigg| >     N \bigg) +1      \,\ge\, c \cdot   \sigma^2   \,.
\end{equation} 
\end{lemma}

\noindent
{\it Proof:} {\bf (I)} \,We argue first \eqref{8.1}, which pertains to the symmetric case. By analogy with \cite{Erd}, p.\,289, and using the symmetry assumption, we obtain for fixed $N \in \N$ the inequalities
$$
 \P \bigg(\,   \sum_{n=1}^N f_n   >   2\, N \bigg) \, \ge \,\P \,\bigg[ \bigcup_{n=1}^N  \bigg( \big\{ f_n > 2\, N \big\} \cap \bigg\{ \sum_{k=1 \atop k \neq n}^N f_k \ge 0 \bigg\} \bigg) \bigg]~~~~~~~~~~~~~~~~
$$
\begin{equation} 
\label{8.3}
~~~~~~~~ ~~~~~~~~\ge \, \sum_{n=1}^N \, \bigg[ \, \P \bigg(   f_n > 2\,N\,,~ \sum_{k=1 \atop k \neq n}^N f_k \ge 0   \bigg) - \P \bigg( \, \bigcup_{k=1 \atop k < n}^N   \big\{ f_n > 2\,N \big\} \cap \big\{ f_k > 2\,N \big\} \bigg) \, \bigg]
\end{equation} 
$$
 \ge \, \sum_{n=1}^N \,  \bigg( \, \frac{1}{\,2\,} \, \P \big( f > 2\,N \big) - N \cdot \Big( \P \big( f > 2\,N \big) \Big)^2 \, \bigg)^+ \,. ~~~~~~~~~~~~~~~~~
$$

\noindent
We distinguish at this point two cases: 
\\ {\bf (a)} If $\, \E  (  | f  |  ) < \infty\,,$ we have $\,\lim_{N \to \infty} \big( N \cdot \P  (  | f  | > N  ) \big) =0\,$ and thus  $\, \P \big( f > 2\,N \big) \le 1 / (4 N  )\,$ for all $ N \ge N_0$ with $N_0$ sufficiently large; then \eqref{8.3} gives $\, \P \big(\,   \sum_{n=1}^N f_n  >   2\, N \big) \, \ge \,(N/4) \cdot \P \big( f > 2\,N \big)\,, $ thus also, for some universal real constant $K >0$, the bound 
$$
\E \big( f^2 \big) \, \le \, K \cdot \sum_{N \in \N}   \frac{N}{4} \cdot \P \big( f > 2\,N \big)     \le K  \cdot \sum_{N \in \N}   \P \Big(\,   \sum_{n=1}^N f_n   >  2\,  N \Big)\,.
$$
{\bf (b)} On the other hand, if $\, \E  (  | f   |  ) = \infty\,,$ the HRE property fails for the sequence $\, f_1, f_2, \cdots\,,$ i.e., we have $\,\E (f^2) = \infty = \P \big(\,   \sum_{n=1}^N f_n  >   2\, N \big)  \,. $

 \smallskip
  In either case,  $\,\sigma^2 \, \le \,       K  \cdot \sum_{N \in \N}   \P \big(\,  \big| \sum_{n=1}^N f_n \big|  >  2\,  N \big)\, $ holds for $\, \sigma^2 = \E \big( f^2 \big) \ge 1\,$ and some universal real constant $K>0$\,; so   \eqref{8.1} holds with   $\,c= 1 / K\,$ after noting that, for $\, \sigma^2   < 1\,,$  it is  satisfied trivially.  We  observe also that this proof gives,   possibly with a different universal constant $c>0$, the inequality                  
 \begin{equation} 
 \label{8.6}
 \sum_{N \in \N      } \,\P \Big(\,  \sum_{n=1}^N f_n   > 2\,N   \Big)\,\ge\,  c   \cdot   \sigma^2\,. 
 \end{equation}    
   
\noindent                
{\bf (II)\,} For the general, non-symmetric case, we argue as follows. We let $\, \big( f_n^\pm \big)_{n \in \N}\,$ be independent copies of the sequence $\, \big( f_n  \big)_{n \in \N}\,$ and set $\, g_n := f^+_n -  f^-_n\,,~ n \in \N\,;$   these $\, \big( g_n  \big)_{n \in \N}\,$ are independent, identically distributed with variance $\, 2 \, \sigma^2 = 2\, \E (f^2)\,$, and {\it symmetric,}  so we obtain from \eqref{8.6} the bound
\begin{equation} 
\label{8.7}
 \sum_{N \in \N  
 } \,\P \Big(\,  \sum_{n=1}^N g_n   > 2\,N   \Big)\,\ge\,  c  \cdot   \E \big( g_1^2 \big) .
\end{equation} 
As   in \cite{Erd}, we note
$\,
\big\{  \sum_{n=1}^N g_n   > 2\,N   \big\}  \, \subseteq \, \big\{\,\big|   \sum_{n=1}^N f^+_n  \big| > N   \big\} \cup \big\{\, \big|  \sum_{n=1}^N f^-_n  \big| > N   \big\}\,,
$ 
therefore  also 
$\, 
\P \big(\,  \sum_{n=1}^N g_n   >  2\,N   \big)  \, \leq \, 2\cdot \P \big(\,\big|   \sum_{n=1}^N f_n  \big| > N  \,\big)\,;
 $
consequently, for some new universal constant $c>0\,,  $ we have  in conjunction with \eqref{8.6}  also   the bound
$$
~~~~~~~~~~~~~
 \sum_{N \in \N  
 } \,\P \Big(\,\Big|   \sum_{n=1}^N f_n  \Big| > N   \Big)\,\ge\, \big( 1 / 2 \big) \, \cdot  \sum_{N \in \N  
 } \,\P \Big(\,  \sum_{n=1}^N g_n   > 2\,N   \Big)\,\ge\, c \, \,  \E \big( f^2 \big) \,. \qquad ~~~~~~~~~~~~~~ \qed
$$

\begin{acks}[Acknowledgments]
Work on this project started in June 2022 during a  visit at the Centro di Ricerca Matematica Ennio \textsc{de\,Giorgi}  of the Scuola Normale Superiore, Pisa. The authors are greatly  indebted to Tomoyuki \textsc{Ichiba} and George \textsc{Stoica} for bringing      important  results to their attention, and  to Amol \textsc{Aggarwal} for his suggestions and his interest in this work.  They would like to thank the anonymous referees, an Associate Editor, and the Editor, for   constructive comments that improved very substantially the quality of the paper.
\end{acks}
 
\begin{funding}
The second author acknowledges support from the National Science Foundation   under      Grants DMS-20-04977 and  DMS-25-06199, as well as  from a Lenfest Award at Columbia University. The third author acknowledges support from the  Austrian Science Fund (FWF) under
 grants P-35519 and P-35197.
\end{funding}

\end{document}